\documentclass[twoside]{article}
\usepackage[accepted]{aistats2021}

\usepackage[round]{natbib}

\usepackage[utf8]{inputenc} % allow utf-8 input
\usepackage[T1]{fontenc}    % use 8-bit T1 fonts
\usepackage{hyperref}       % hyperlinks
\usepackage{url}            % simple URL typesetting
\usepackage{booktabs}       % professional-quality tables
\usepackage{amsfonts}       % blackboard math symbols
\usepackage{nicefrac}       % compact symbols for 1/2, etc.
\usepackage{microtype}
\usepackage{amssymb, amsmath, amsthm}
\usepackage{microtype}

\usepackage{graphicx}
\usepackage{subcaption}
\usepackage{booktabs} % for professional tables
\usepackage{algorithm}
\usepackage{algorithmic}
\usepackage{color}
\usepackage{tcolorbox}
\tcbset{boxsep=0mm,boxrule=0pt,colframe=white,arc=0mm,left=0.5mm,right=0.5mm}
\usepackage{thmtools}
\usepackage{thm-restate}
\usepackage{enumerate}
\usepackage{enumitem}
\usepackage{url}
\usepackage{wrapfig}
\usepackage{dirtytalk}

\DeclareMathOperator*{\argmin}{argmin}
\newcommand{\grad}{\text{grad}}
\newcommand{\Hess}{\text{Hess}}

\renewcommand{\H}{{\mathcal H}}

\newcommand{\R}{{\mathbb{R}}}

\newtheorem{theorem}{Theorem}
\newtheorem{lemma}[theorem]{Lemma}
\newtheorem{proposition}[theorem]{Proposition}

\newtheorem{assumption}{Assumption}
\newtheorem{definition}{Definition}

\newcommand{\bigO}{\mathcal{O}}

\begin{document}

\twocolumn[

\aistatstitle{Momentum Improves Optimization on Riemannian Manifolds}

\aistatsauthor{ Foivos Alimisis \And Antonio Orvieto \And  Gary B{\'e}cigneul \And Aurelien Lucchi }

\aistatsaddress{ IST Austria \And  ETH Z{\"u}rich \And Gematria Technologies \\
London, U.K. \And ETH Z{\"u}rich } ]

\begin{abstract}
We develop a new Riemannian descent algorithm that relies on momentum to improve over existing first-order methods for geodesically convex optimization. In contrast, accelerated convergence rates proved in prior work have only been shown to hold for geodesically strongly-convex objective functions. We further extend our algorithm to geodesically weakly-quasi-convex objectives. Our proofs of convergence rely on a novel estimate sequence that illustrates the dependency of the convergence rate on the curvature of the manifold. We validate our theoretical results empirically on several optimization problems defined on the sphere and on the manifold of positive definite matrices.
\end{abstract}

\section{Introduction}
The field of optimization plays a central role in machine learning. At its core lies the problem of finding a minimum of a function $f:H \rightarrow \R$. In the vast majority of applications in machine learning, $H$ is considered to be a Euclidean vector space. However, a number of machine learning tasks can profit from a specialized problem-dependent Riemannian structure \citep{bonnabel2013stochastic, zhang2016first}, which will be the focus of our discussion in this paper. 
Among the most popular types of methods to optimize $f$ are first-order methods, such as gradient descent that simply updates a sequence of iterates $\{ x_k \}$ by stepping in the opposite direction of the gradient $\nabla f(x_k)$. In the case $H = \mathbb{R}^n$, gradient descent as a first-order method has been shown to achieve a suboptimal convergence rate on convex problems ($\bigO(1/k )$). In a seminal paper, \citep{nesterov1983method}, Nesterov showed that one can construct an optimal --- i.e. \textit{accelerated} --- algorithm that achieves faster rates of convergence for both convex ($\bigO(1/k^2)$) and strongly-convex functions. The convergence analysis of this algorithm relies heavily on the linear structure of $H$ and it is not until recently that a first adaptation to Riemannian manifolds was derived by \cite{zhang2018towards}. The algorithm by \cite{zhang2018towards} is shown to obtain an accelerated rate of convergence for functions that are known to be \emph{geodesically} strongly-convex, provided that one initializes in a neighborhood of the (unique) solution. These functions are of particular interest as they might be non-convex in the Euclidean sense and they occur in some relevant computational tasks, such as the approximation of the Karcher mean of positive definite matrices~\citep{zhang2016riemannian}. However, many other interesting problems belong to the weaker class of \emph{geodesically} convex functions. This includes problems defined on the cone of Hermitian positive definite matrices~\citep{sra2015conic}, which appear in various areas of machine learning such as tracking~\citep{cheng2013novel} and medical imaging~\citep{zhu2007statistical}.
In this paper, we therefore address the question of whether an algorithm that relies on momentum can provably achieve a faster rate of convergence for functions that are geodesically convex but not necessarily strongly convex. We also consider the extension to the \emph{weaker} class of geodesically weakly-quasi-convex objective functions. A more thorough motivation for investigating convex and weakly-quasi convex objectives in Riemannian optimization can be found in Section 4 of~\citep{alimisis2019continuoustime}. Our main contributions are:
\begin{enumerate}[wide, labelwidth=1pt, labelindent=1pt]
    \itemsep0em
    \item We propose a new \emph{practical} Riemannian algorithm that exploits momentum to speed up convergence for geodesically convex and weakly-quasi-convex functions.
    As in~\citep{nesterov2018primal}, our approach uses a small-dimensional relaxation~(sDR) oracle (which can be solved \textit{approximately} and in linear time) to perform adaptive linear coupling\footnote{See discussion in the next section.}~\citep{allen2014linear}.  In order to provide theoretical guarantees for this new algorithm, we use a novel estimate sequence combining ideas from \citep{nesterov2018primal} and \citep{zhang2018towards} as well as develop some new results at the intersection of optimization and Riemannian geometry.
    \itemsep0em
    \item Our main algorithm applied to geodesically convex objective functions provides better theoretical guarantees of convergence compared to Riemannian Gradient Descent (RGD)~\citep{zhang2016first}, given that the bound on the working domain is not exceedingly large. Since RGD is the only known first-order method with guaranteed convergence for geodesically convex functions, our algorithm has the best known worst-case behaviour. Moreover, our algorithm is accelerated for the first (practically large) part of the optimization procedure. 
    \itemsep0em
    \item We validate our theoretical findings numerically on several important machine learning problems defined on manifolds of both positive curvature~(Rayleigh quotient maximization) and negative curvature~(operator scaling and Karcher mean approximation). Some of these problems are convex (but not strongly-convex) while others have a relatively small strong convexity constant. We show the empirical superiority of our method when compared to Riemannian algorithms designed for well-conditioned geodesically strongly-convex objectives, such as RAGD~\citep{zhang2018towards}.
\end{enumerate}
\section{Related Work}
\paragraph{Accelerated Gradient Descent (AGD).}
The first accelerated gradient descent algorithm in Euclidean vector spaces is due to~\cite{nesterov1983method}. Since then, the community has shown a deep interest in understanding the mechanism underlying acceleration. A recent trend has been to look at acceleration from a continuous-time viewpoint~\citep{su2014differential, wibisono2016variational}. In this framework, AGD is seen as the discretization of a second-order ODE.
Alternatively, \cite{allen2014linear} showed how one can view AGD as a primal-dual method performing linear coupling between gradient descent and mirror descent. Recently~\cite{nesterov2018primal} proposed AGDsDR, a modification of the method by~\cite{allen2014linear}, which adaptively selects the linear coupling parameter~(denoted by $\beta$) at each iteration using an approximate line search. This work will serve as an inspiration for us to design an accelerated Riemannian algorithm.
\vspace{-2mm}
\paragraph{Riemannian optimization.} Research in the field of Riemannian optimization has encountered a lot of interest in the last decade. A seminal book in the area is~\citep{absil2009optimization}, which gives a comprehensive review of many standard optimization methods, but does not discuss acceleration. More recently, \cite{zhang2016first} proved convergence rates for Riemannian gradient descent applied to geodesically convex functions. Acceleration in a Riemannian framework was first discussed by~\cite{liu2017accelerated}, who claimed to have designed a Riemannian method with guaranteed acceleration. While their methodology is interesting, unfortunately, as discussed in~\citep{zhang2018towards}, their algorithm relies on finding the \textit{exact} solution to a nonlinear equation at each iteration, and it is not clear how difficult this additional problem might be or how approximation errors accumulate. Subsequently, \cite{zhang2018towards} developed the first computationally tractable accelerated algorithm on a Riemannian manifold, but their approach only has provable convergence for geodesically strongly-convex objectives~(provided that one initializes sufficiently close to the solution). A more recent work, \citep{ahn2020nesterov}, attempts to tackle the problem of acceleration for geodesically strongly-convex optimization with global convergence rate~(no assumptions on the initialization required). Notwithstanding the theoretical significance of this work, the final algorithm has the practical drawback that achieves full acceleration only in late training~(after a possibly very large number of steps for ill-conditioned problems), while at the beginning behaves comparably to Riemannian gradient descent. Instead, using a continuous-time viewpoint, the recent work \citep{alimisis2019continuoustime} analyzed various ODEs that can model acceleration on Riemannian manifolds with theoretical guarantees of convergence. They derived discrete-time algorithms via numerical integration of the continuous-time process but do not provide theoretical guarantees for the discrete-time schemes. 
The problem we address is different from prior work as we aim to demonstrate that momentum provably yields a better rate of convergence than Riemannian gradient descent for the classes of geodesically convex and weakly-quasi-convex functions, which are both of significant practical interest (see discussion in Section~\ref{sec:experiments}). We note that extending the proof by~\cite{zhang2018towards} to these weaker classes of functions is not straightforward due to some distortions between the tangent spaces of the sequence of iterates of the algorithm~\footnote{By "distortion", we mean that when considering two successive iterates $x_k$ and $x_{k+1}$, the terms $\log_{x_k}(a)-\log_{x_k}(b)$ and $\log_{x_{k+1}}(a)-\log_{x_{k+1}}(b)$ appearing in the estimate sequence belong to different tangent spaces and are therefore not directly comparable (while they are exactly the same in the Euclidean case).}. Indeed, the estimate sequence used in \citep{zhang2018towards} relies on changing the tangent space at each step. These successive changes give rise to additional errors which can be dealt with by relying on the strong convexity assumption. However, we were unable to adapt their proof to weaker function classes. Instead, we rely on a novel estimate sequence that is qualitatively different from the one used in \citep{zhang2018towards} in order to avoid distortions produced by changing tangent spaces.

\section{Background}
\subsection{Preliminaries from Differential Geometry}
We review some basic notions from Riemannian geometry that are required in our analysis. For a full review, we refer the reader to some classical textbook, for instance \citep{spivak}.
\vspace{-2mm}
\paragraph{Manifolds.}
A differentiable manifold $M$ is a topological space that is locally Euclidean. This means that for any point $x \in M$, we can find a neighborhood that is diffeomorphic to an open subset of some Euclidean space. This Euclidean space can be proved to have the same dimension, regardless of the chosen point, called the dimension of the manifold.
\newline
A Riemannian manifold $(M,g)$ is a differentiable manifold equipped with a Riemannian metric $g_x$, i.e.  an inner product for each tangent space $T_xM$ at $x \in M$. We denote the inner product of $u,v \in T_x M$ with $\langle u,v \rangle_x$ or just $\langle u,v \rangle$ when the tangent space is obvious from context. Similarly we consider the norm as the one induced by the inner product at each tangent space. The Riemannian metric provides us a way to measure the distance $d$ between points on the manifold, transforming it into a metric space. Given $A \subseteq M$, the diameter of $A$ is defined as $\textnormal{diam}(A)=\sup_{p,q \in A} d(p,q)$.
\vspace{-5mm}
\paragraph{Geodesics.}
Geodesics are curves $\gamma: [0,1] \rightarrow M$ of constant speed and of (locally) minimum length. They can be thought of as the Riemannian generalization of straight lines in Euclidean space. Geodesics are used to construct the exponential map $ \exp_x: T_x M \rightarrow M$, defined by $ \exp_x(v)= \gamma(1)$, where $\gamma$ is the unique geodesic such that $\gamma(0)=x$ and $\dot \gamma(0)=v$. The exponential map is locally a diffeomorphism. We denote the inverse of the exponential map $\exp_x$ (in a neighborhood $U \subseteq M$ of $x$) by $\log_x: U \to T_x M$. Geodesics also provide a way to transport vectors from one tangent space to another. This operation, called parallel transport, is usually denoted by $\Gamma_x^y: T_x M \rightarrow T_y M$.
\vspace{-2mm}
\paragraph{Vector fields and covariant derivative.}
The correct notion to capture second order changes on a Riemannian manifold is called covariant differentiation and it is induced by the fundamental property of Riemannian manifolds to be equipped with a connection. We are interested in a specific type of connection, called the Levi-Civita connection, which induces a specific type of covariant derivative. The fact that a unique Levi-Civita connection exists always in a Riemannian manifold is the subject of the fundamental theorem of Riemannian geometry. However, for the purpose of our analysis, it will be sufficient to rely on a simple notion of covariant derivative that relies on the (more visualizable) notion of parallel transport. First, we define vector fields on a Riemannian manifold as sections of the tangent bundle.

\begin{definition}
Let $M$ be a Riemannian manifold. A vector field $X$ in $M$ is a smooth map $X:M \rightarrow \mathcal{T}M$, where $\mathcal{T}M$ is the tangent bundle, i.e. the collection of all tangent vectors in all tangent spaces of $M$, such that $p \circ X$ is the identity ($p$ projects from $\mathcal{T}M$ to $M$).
\end{definition}
\vspace{-2mm}

One can see a vector field as an infinite collection of imaginary curves, the so-called integral curves (formally solutions of first order differential equations on $M$).
\begin{definition}
Given two vector fields $X,Y$ in a Riemannian manifold $M$, we define the covariant derivative of $Y$ along $X$ to be
$\nabla_X Y (p) = \lim_{h\to0}\frac{\Gamma_{\gamma(h)}^{p} Y(\gamma(h))-Y(p)}{h}$, where $\gamma$ is the unique integral curve of $X$, starting from $p$, i.e $\gamma(0)=p$.
\label{def:covariant_derivative}
\end{definition}

\subsection{Geodesic convexity}
We remind the reader of the basic definitions needed in Riemannian optimization.
\begin{definition}
A subset $A\subseteq M$ of a Riemannian manifold $M$ is called geodesically uniquely convex, if every two points in $A$ are connected by a unique geodesic.
\end{definition}
\begin{definition}
A function $f:A \rightarrow \R$ is called geodesically convex, if for any $p,q \in M$, we have $f(\gamma(t)) \leq (1-t)f(p)+tf(q)$ for any $t \in [0,1]$, where $\gamma$ is the geodesic connecting $p,q \in M$. 
\end{definition}

Given a function $f:M \rightarrow \R$, the notions of  differential and (Riemannian) inner product allow us to define the Riemannian gradient of $f$ at $x \in M$, which is a tangent vector belonging to the tangent space based at $x$, $T_x M$.
\vspace{-2mm}
\begin{definition}
The Riemannian gradient $\textnormal{\textnormal{gradf}}$ of a (real-valued) function $f:M \rightarrow \R$ at a point $x \in M$, is the tangent vector at $x$, such that $\langle \textnormal{\textnormal{gradf}}(x),u \rangle = df(x)u$~\footnote{$df$ denotes the differential of $f$, i.e. $df (x)[u] = \lim_{t \to 0} \frac{f(c(t)) - f(x)}{t},$ where $c: I \to M$ is a smooth curve such that $c(0) = x$ and $\dot c(0) = u$.}, for any $u \in T_x M$.
\end{definition}
Given the notion of Riemannian gradient and covariant derivative we define the notion of Riemannian Hessian.
\vspace{-2mm}
\begin{definition}
The Hessian of $f$ is defined as a bilinear form at each point $p \in M$, given by $\operatorname{Hess}_p(f)(X,Y)= \langle \nabla_X \operatorname{grad}f,Y \ \rangle$,  for two vector fields $X,Y$ on $M$.
\end{definition}

Using the Riemannian inner product and gradient, we can formulate an equivalent definition for geodesic convexity for a smooth function $f$ defined in a geodesically uniquely convex domain $A$.
\begin{proposition}
Let $f:A \rightarrow \R$ be a smooth, geodesically convex function. Then, for any $x,y \in A$,
\begin{align*}
    f(y)-f(x) \geq \langle \textnormal{\textnormal{gradf}}(x), \textnormal{\textnormal{log}}_x(y) \rangle.
\end{align*}
\end{proposition}

As in the Euclidean case, any local minimum of a geodesically convex function is a global minimum.
We now generalize the well-known notion of Euclidean weak-quasi-convexity to Riemannian manifolds. For a review of this notion, we refer the reader to \citep{guminov2017accelerated}.
\begin{definition}A function $f:A \rightarrow \R$ is called geodesically $\alpha$-weakly-quasi-convex with respect to $c \in M$, if $\alpha (f(c)-f(x)) \geq \langle \textnormal{gradf}(x), \textnormal{log}_x(c) \rangle$ for some fixed $\alpha \in (0,1]$ and any $x \in A$.
\end{definition}

It is easy to see that weak-quasi-convexity implies that any local minimum of $f$ is also a global minimum. Using the notion of parallel transport we can define when $f$ is geodesically $L$-smooth, i.e. has Lipschitz continuous gradient in a differential-geometric way.
\begin{definition} $f:A \rightarrow \R$ is called L-smooth if $\| \textnormal{gradf}(x) -\Gamma_y^x {\textnormal{gradf}(y)} \| \leq L \| \textnormal{log}_x(y) \|$ for any $x,y \in A$.
\end{definition}

Geodesic $L$-smoothness has similar properties to its Euclidean analogue: a two times differentiable function is $L$-smooth, if and only if the norm of its Riemannian Hessian is upper bounded by $L$.

\subsection{Basic Assumptions}

In this paper, we make the standard assumption that the input space is not "infinitely curved". In order to make this statement rigorous, we need the notion of \textit{sectional curvature} $K$, which is a measure of how sharply the manifold is curved (or how \say{far} from being flat our manifold is), \say{two-dimensionally}. More concretely, as in~\citep{zhang2018towards}, we make the following set of assumptions:
\begin{assumption}
Given $A \subseteq M$ geodesically uniquely convex, and $f: A \to \R$,
\begin{itemize}[leftmargin=.25in]
    \setlength\itemsep{0mm}
    \item[1.] The sectional curvature $K$ inside $A$ is bounded from above and below, i.e. $ K_{min} \leq K \leq K_{\max}$.
    \item[2.] $A$ is a geodesically uniquely convex subset of $M$, such that $\textnormal{diam}(A) \leq D < \infty$. This implies that the exponential map $\exp_x:T_x M \rightarrow M$ is globally a diffeomorphism for any $x \in A$ with inverse denoted by $\log_x$. 
    \item[3.] $f$ is geodesically $L$-smooth with its local minima (which are all global) inside $A$, we denote some of them by $x^*$.
    \item[4.] We have granted access to oracles which compute the exponential and logarithmic maps as well as the Riemannian gradient of $f$ efficiently.
    \item[5.] All the iterates of our algorithms remain inside $A$.
\end{itemize}
\label{ass:main_assumption}
\end{assumption}

The last assumption is standard in accelerated Riemannian optimization, \citep{zhang2018towards,alimisis2019continuoustime,ahn2020nesterov}, and we did not observe it to be violated in our experiments. However, it remains an open question as to whether it could be relaxed or even removed completely from our analysis.

%%%%%%%%%%%%%%%%%%%%%%%%%%%%%%%%%%%%%%%%%%%%%%%%%%%%%%%%%%%%
%%%%%%%%%%%%%%%%%%%%%%%%%%%%%%%%%%%%%%%%%%%%%%%%%%%%%%%%%%%%

\section{The RAGDsDR Algorithm}
\begin{figure}
\centering
        \includegraphics[width=0.85\linewidth]{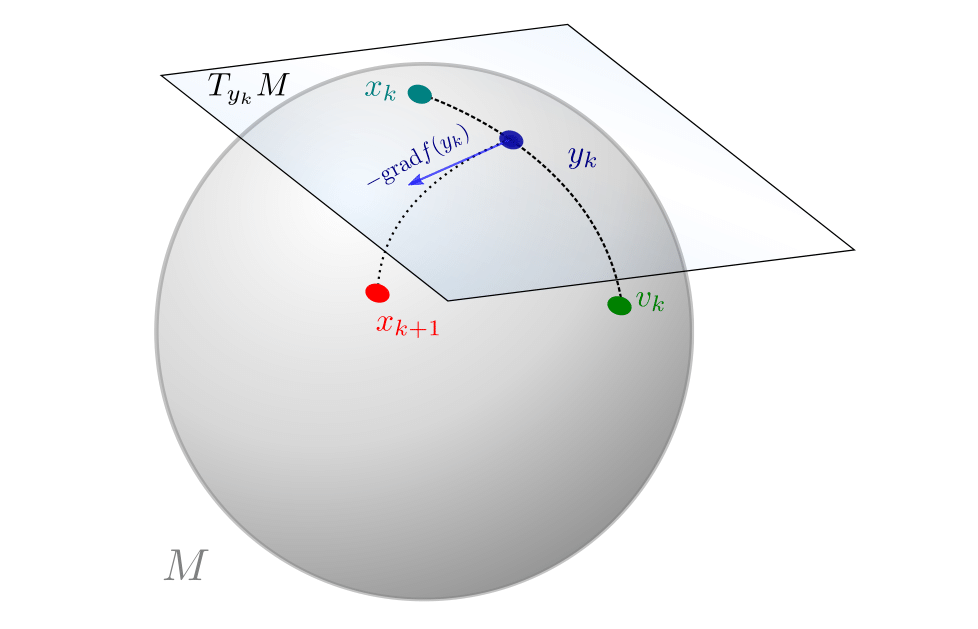}
        \caption{\footnotesize{Illustration of one step of Algorithm \ref{alg:RAGDsDR}. The point $y_k$ is computed to minimize $f$ on the geodesic between $x_k$ and $v_k$. $x_{k+1}$ is computed by taking a gradient step from $y_k$ and projected using the exponential map $\exp_{y_k}(\cdot)$.}}
         \label{fig:algorithm}
\end{figure}

We now develop a new Riemannian algorithm that relies on momentum and which is inspired by the Euclidean algorithm presented in \citep{nesterov2018primal} (see description in Appendix \ref{app:appendix A}). It is detailed in Algorithm \ref{alg:RAGDsDR} and illustrated in Figure \ref{fig:algorithm}. At each iteration $k$, the next iterate $x_{k+1}$ (line 5) is computed by taking a gradient step at an interpolated point $y_k$ (line 4) which follows the direction of a momentum term $\log_{v_k} (x_k)$. The main difference with the Euclidean case is that the curve from $v_k$ to $x_k$ is a geodesic on the manifold $M$ instead of a straight Euclidean line. As in~\citep{nesterov2018primal}, we also rely on a minimization over a closed interval (i.e. the small-dimension relaxation, sDR) to choose the best possible stepsize $\beta_k$ (line 3) on the geodesic connecting $v_k$ to $x_k$. We will see in the next section that this minimization is computationally \textit{fast to solve}~(also see Section \ref{sec:experiments}), can be computed approximately and practically yields faster convergence than the typical fixed parameter $\beta_k = \frac{k}{k+2}$~\citep{nesterov2018lectures}. The curvature of $M$ is involved directly in the algorithm via the quantity $\zeta \geq 1$ (line 6), defined as 
\begin{equation}
\label{eq:zeta}
    \zeta:= 
    \begin{cases}
    \sqrt{-K_{\min}} D \coth(\sqrt{-K_{\min}} D) &, K_{\min}<0 \\
    1 &, K_{\min}\geq 0
    \end{cases}
\end{equation}

\begin{algorithm}[H]
\small
\caption{RAGDsDR for convex functions}
\begin{algorithmic}[1]
\label{alg:RAGDsDR}
      \STATE $A_0=0,x_0=v_0 \in A$
      \FOR {$k \geq 0$}
      \STATE $\beta_k =  \underset{\beta \in [0,1]}{\mathrm{argmin}} \left \lbrace f(\exp_{v_k} (\beta \log_{v_k}(x_k))) \right \rbrace$
      \STATE $y_k=\exp_{v_k}(\beta_k \log_{v_k}(x_k))$
      \STATE $x_{k+1} = \exp_{y_k}\left(-\frac{1}{L} \grad f(y_k)\right)$
      \STATE $\frac{\zeta a_{k+1}^2}{A_k+a_{k+1}}=\frac{1}{L}, a_{k+1}>0$
      \STATE $A_{k+1}=A_k+a_{k+1}$
      \STATE $v_{k+1}=\exp_{v_k} \left(-a_{k+1} \Gamma_{y_k}^{v_k} \grad f(y_k)\right)$ 
      \ENDFOR
\end{algorithmic}      
 \end{algorithm}

The discriminant of the quadratic equation defining $a_{k+1}$ at step 6 is positive, thus the aforementioned equation has a positive and negative solution, from which we choose the first.

The computation of the parallel transport at step 8 is given directly by the oracle, since it relies on the computation of the exponential map. For manifolds found in applications, the parallel transport is cheap and implementations are found in libraries such as~\citep{townsend2016pymanopt}. 

The definition of $v_{k+1}$ at step 8 is qualitatively different from the one in \citep{zhang2018towards}, but not heavier computationally, since in both cases we need three oracle calls. Note that~\cite{nesterov2018primal} define $v_{k+1}$ through a minimization problem (see Appendix~\ref{app:appendix A}). 
This approach could be naively generalized to the Riemannian setting but it would yield a minimization problem that has no explicit solution due to non-linearity. Instead, we find a generalization of the Euclidean definition that can be solved explicitly and write $v_{k+1}$ directly in its explicit form in Algorithm \ref{alg:RAGDsDR}.

\paragraph{Geodesic search.} The second step in Algorithm \ref{alg:RAGDsDR} is solved using a procedure similar to a line search which we name \textit{geodesic search}. It guarantees that the following two key conditions hold (proof in App.~\ref{app:appendix line search}):
\begin{align}
f(y_k) \leq f(x_k) \text{ and } \langle \grad f(y_k), \log_{y_k}(v_k) \rangle \geq 0.
\label{eq:line_search_cond}
\end{align}

Practically, the geodesic search procedure is inexact. While we can still expect the first inequality in Eq. \ref{eq:line_search_cond} to be satisfied exactly, the second one can only be satisfied up to a small error $\Tilde{\epsilon}>0$, i.e. $\langle \grad f(y_k), \log_{y_k}(v_k) \rangle \geq -\Tilde{\epsilon}$. We note that this is an analogous condition to the one used by~\cite{nesterov2018primal} in the Euclidean case. As we will see shortly, one of the main quantities of interest in our analysis will be
\begin{equation}
\mathcal{E}_k(x) := \langle \grad f(y_k), \log_{y_k}(x)- \Gamma_{v_k}^{y_k} \log_{v_k}(x) \rangle,
\label{eq:error_term}
\end{equation}
which occurs as an error in our estimate-sequence analysis and captures the curved nature of the manifold $M$. We will prove that the absolute value of this error is bounded by the sum of two terms, namely $\mid \mathcal{E}_k(x) \mid \leq \Tilde{\epsilon}+\Tilde{\eta}_k$, where $\Tilde{\epsilon}$ is the error obtained by the geodesic search and $\Tilde{\eta}_k$ is an extra curvature-dependent error. The latter depends on an upper bound on the working domain $D$ and it decays to $0$ as the algorithm runs. In the Euclidean case $\Tilde \eta_k=0$. We will prove that in the Riemannian case $\Tilde{\eta}_k=\bigO\left(\frac{d(M)}{k}\right)$, where $d(M)$ is a small constant which depends on the sectional curvature and a bound of our working domain.

%%%%%%%%%%%%%%%%%%%%%%%%%%%%%%%%%%%%%%%%%%%%%%%%%%%%%%%%%%%%
%%%%%%%%%%%%%%%%%%%%%%%%%%%%%%%%%%%%%%%%%%%%%%%%%%%%%%%%%%%%

\section{Convergence Analysis}

\paragraph{Geodesically-convex functions.} We now present our main convergence result.

Our analysis is based on a novel estimate sequence, which allows for an extra error at each step. However, this extra error does not accumulate, and it decays linearly over iterations. As a result, we obtain a rate of convergence that is superior to the convergence guarantees of RGD derived in \citep{zhang2016first} under a restriction on the bound of the working domain (see later discussion).
\newline
We first need to examine the behaviour of $\mathcal{E}_k(x)$:
\begin{lemma}
\label{le:extra_error}
Under our set of assumptions (Assumption \ref{ass:main_assumption}), Algorithm~\ref{alg:RAGDsDR} produces iterates $y_k, v_k$ such that
\begin{equation*}
    -\mathcal{E}_k(x) \leq \| \grad f(y_k) \| \max \lbrace \zeta-1,1-\delta \rbrace D+\Tilde{\epsilon}
\end{equation*}
with $\zeta \geq 1$ defined by equation (\ref{eq:zeta}) and $\delta \leq 1$ defined by
\begin{equation*}
    \delta:= 
   \begin{cases}
    1 &, K_{\max} \leq 0\\
    \sqrt {K_{\max}} D \cot(\sqrt{K_{\max}} D) &, K_{\max} > 0 
   \end{cases}
\end{equation*}
\end{lemma}

We prove this lemma in Appendix \ref{app:extra_error}. We rely on various geometric bounds derived in Appendix \ref{app:operator}, which are inspired by those of \cite{alimisis2019continuoustime}. Generally speaking, $\delta$ and $\zeta$ are obtained by considering the spectral properties of an operator similar to the Riemannian Hessian of the squared distance function, $\delta$ as a lower bound of its smallest eigenvalue and $\zeta$ as an upper bound of the largest one.

We are now ready to state our main convergence result:
\begin{tcolorbox}
\begin{theorem}
\label{thm:main_theorem}
Algorithm \ref{alg:RAGDsDR} applied to a geodesically convex function $f$ produces iterates $x_k$, such that
\begin{align*}
    &f(x_k)-f^*  \leq \\ & \frac{2 \zeta L D^2}{k^2}+4 \max \lbrace\zeta-1,1-\delta\rbrace \frac{ \zeta L D^2}{k}+\Tilde{\epsilon}   \leq \\ & 2 \max\left \lbrace \frac{2 \zeta L D^2}{k^2},4 \max  \lbrace\zeta-1,1-\delta\rbrace \frac{ \zeta L D^2}{k} \right \rbrace+\Tilde{\epsilon}
\end{align*}
\end{theorem}
\end{tcolorbox}
Recall that parameter $D$ is used to denote an upper bound for the diameter of our working domain (Assumption \ref{ass:main_assumption}).

The proof is derived in Appendix \ref{app:Main Analysis} and relies on Lemma \ref{le:extra_error}. At first glance, the upper bound seems rather intuitive for those familiar with Riemannian optimization, namely the positive-curvature case provides the same guarantees as the Euclidean one, while the negative-curvature case provides worse guarantees.
\newline
Let's take a closer look at the rate of convergence. To do so, we define the following quantity which we call the \say{discrepancy} of the manifold $M$.
\begin{definition}
The discrepancy of the manifold $M$ is defined as $d(M):= 4 \max  \lbrace\zeta-1,1-\delta\rbrace$.
\end{definition}
In the Euclidean case, we have $d(M)=0$, thus our algorithm is a generalization of accelerated gradient descent with line search.
\newline
The convergence rate in Thm.~\ref{thm:main_theorem} is accelerated when
\begin{equation*}
    \frac{2 \zeta L D^2}{k^2} \geq d(M) \frac{ \zeta L D^2}{k}
\end{equation*}
which is equivalent to $k \leq \cfrac{2}{d(M)}.$
Thus, $2/d(M)$ is an upper bound indicating how many steps of the algorithm can be performed with an accelerated convergence rate. When the manifold $M$ tends to be Euclidean in the sense that
$\max\lbrace \mid K_{min} \mid , \mid K_{max} \mid \rbrace \rightarrow 0$, then $d(M) \rightarrow 0$ and $\frac{2}{d(M)} \rightarrow \infty$, increasing the numbers of iterations that one can perform accelerated optimization.
\newline
Even when we exceed this bound, the condition 
\begin{equation}
\label{cond:discrepancy}
    2 d(M) < \frac{1}{2} \Leftrightarrow \max \lbrace\zeta-1,1-\delta\rbrace < \frac{1}{16} 
\end{equation}
suffices to guarantee a better worst-case upper bound than Riemannian Gradient Descent in \citep{zhang2016first} (Theorem 13), since we have a smaller constant in the numerator. This is because the rate provided in Theorem 13 in \citep{zhang2016first} is
\begin{equation*}
    f(x_k)-f^* \leq \frac{\zeta L D^2}{2(\zeta+k-2)}
\end{equation*}
Condition \eqref{cond:discrepancy} implies that $\frac{\zeta L D^2}{2(\zeta+k-2)}>2d(M) \frac{ \zeta L D^2}{k}$, since $4 d(M) <1$ and $\zeta<2$ (which implies that \newline $\frac{1}{k} <\frac{1}{\zeta+k-2}$).
\newline

Finally, condition \eqref{cond:discrepancy} is for instance satisfied if
$$\sqrt{\mid K_{min} \mid }D \leq 0.4,\quad \text{and}\quad \sqrt{\mid K_{max} \mid }D \leq 0.4.$$
Both conditions hold if and only if the curvature of the manifold is in absolute value less or equal than $\frac{0.16}{D^2}$. We summarize these facts in the following theorem:
\begin{tcolorbox}
\begin{theorem}
\label{thm:flex}
When the sectional curvature $K$ of the manifold $M$ satisfies
\vspace{-2mm}
\begin{equation*}
    \mid K \mid \leq \frac{0.16}{D^2},
\end{equation*}
Algorithm \ref{alg:RAGDsDR} performs better than RGD in \citep{zhang2016first}. 
\newline
When 
\vspace{-2mm}
\begin{equation*}
    k \leq \frac{2}{d(M)} \underset{K \rightarrow 0}{\longrightarrow} \infty,
\end{equation*}
Algorithm \ref{alg:RAGDsDR} is accelerated.
\end{theorem}
\end{tcolorbox}
In practical situations, we have observed that the quantity $d(M)$ is very small, and we therefore empirically observe acceleration for a very large number of iterations (almost until convergence). We refer the reader to the discussion in Section~\ref{sec:discussion} which also includes a comparison with \citep{zhang2018towards}.

\paragraph{Geodesically weakly-quasi-convex functions.}
We extend Algorithm~\ref{alg:RAGDsDR} to functions that are $\alpha$-weakly-quasi-convex. This requires to restart Algorithm \ref{alg:RAGDsDR} whenever the suboptimality at the current iteration is less than the previous one by a factor $1-\frac{\alpha}{c}$, where $c>1$ is a constant. This procedure yields Alg.~\ref{alg:RAGDsDR_weak}.

\begin{tcolorbox}
\begin{restatable}{theorem}{ConvergenceTheoremQuasiConvex}
\label{th:convergence_thm_quasi_convex}
Algorithm 2 applied to an $\alpha$-weakly-quasi-convex function as in the assumptions produces a sequence of iterates $\{x_k\}_{k=1}^N$, such that
\begin{align*}
    &f(x_N)-f(x^*) \leq \\ &\bigO \left(\frac{\zeta L D^2 }{\alpha^3 N^2} \right)+d(M) \bigO\left(\frac{\zeta L D^2}{\alpha^2 N} \right)  +\frac{c}{(c-1)\alpha}\Tilde{\epsilon},
\end{align*}
where $\Tilde{\epsilon}$ is the error of the geodesic search, $c>1$ and $d(M)$ is the discrepancy of the manifold.
\end{restatable}
\end{tcolorbox}

\begin{algorithm}[h]
\caption{RAGDsDR for weakly-quasi-convex functions}
\label{alg:RAGDsDR_weak}
\begin{algorithmic}[1]
\FOR {$i \geq 0$}
     \STATE $A_0=0,x_0^i=v_0^i \in A$
     \FOR {$k \geq 0$}
     \STATE {$\beta_k =$\small $\argmin_{\beta \in [0,1]} \left \lbrace f \left(\exp_{v_k^i} \left(\beta \log_{v_k^i}(x_k^i)\right)\right) \right \rbrace$}
     \STATE $y_k^i=\exp_{v_k^i}\left(\beta_k \log_{v_k^i}(x_k^i)\right)$
     \STATE $x_{k+1}^i = \exp_{y_k^i}\left(-\frac{1}{L} \grad f(y_k^i) \right)$
     \STATE $\frac{\zeta a_{k+1}^2}{A_k+a_{k+1}}=\frac{1}{L}, a_{k+1}>0$
     \STATE $A_{k+1}=A_k+a_{k+1}$
     \STATE $v_{k+1}^i=\exp_{v_k^i}\left(-a_{k+1} \Gamma_{y_k^i}^{v_k^i} \grad f(y_k^i)\right)$
     \IF {$f(x_k^i)-f(x^*) \leq \left(1-\frac{\alpha}{c}\right) (f(x_0^i)-f(x^*))$}
     \STATE \textbf{break}
     \ENDIF
     \ENDFOR
     \STATE $x_0^{i+1}=x_N^i$ (\textnormal{where $N$ is the number of steps performed in the loop over $k$})
     \ENDFOR
\end{algorithmic}
\end{algorithm}
As in the convex case, the $\bigO(\frac{1}{N})$ part of the rate is multiplied by the discrepancy of the manifold $d(M)$, thus the analysis of Theorem \ref{thm:flex} holds almost the same. The proof can be found in Appendix~\ref{app:weakly-quasi-convex}.

\section{Numerical Experiments}
\label{sec:experiments}
We validate our findings on Riemannian manifolds of both positive and negative curvature. Our code\footnote{\url{https://github.com/aorvieto/RAGDsDR}}
is built on top of \texttt{PyManopt}~\citep{townsend2016pymanopt}. We compare RAGDsDR~(Algorithm~\ref{alg:RAGDsDR}) with Riemannian Gradient Descent~(RGD) and, when possible (i.e. when we can estimate the strong convexity modulus), with RAGD by~\cite{zhang2018towards}. As a more practical alternative to the geodesic search in step 2~(which we solve with at most 10 iterations of \textit{golden-section search}), we show the performance for $\beta_k = \frac{k}{k+2}$. Under this choice, RAGDsDR recovers a Riemannian version of Linear Coupling~\citep{allen2014linear}.

\subsection{Positive curvature}
\label{sec:experiments_positive}

\begin{figure}[H]
\begin{center}
\begin{subfigure}{0.54\linewidth}
  \centering
  \includegraphics[width=\textwidth]{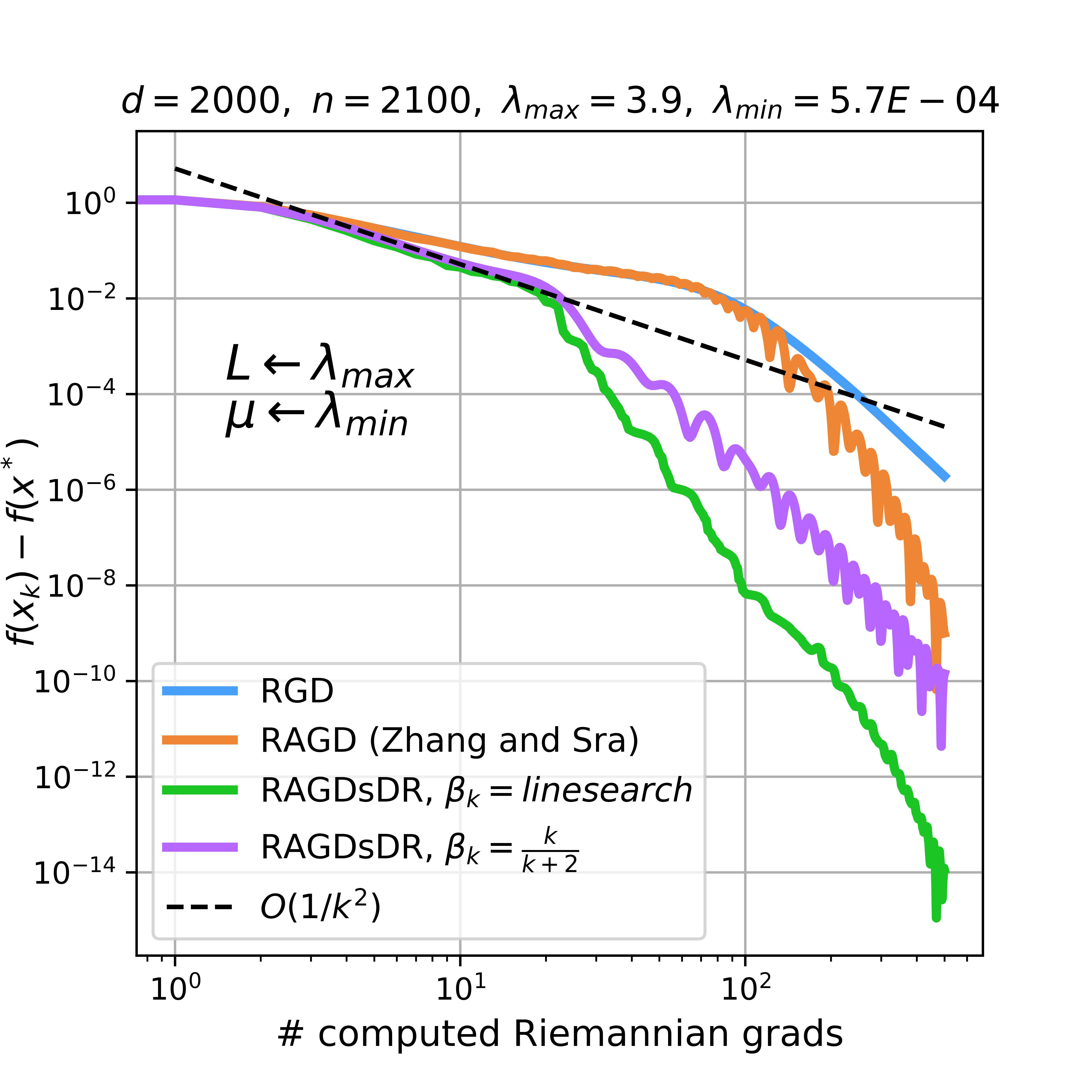}
\end{subfigure}%
\begin{subfigure}{.4\linewidth}
  \centering
  \includegraphics[width=\textwidth]{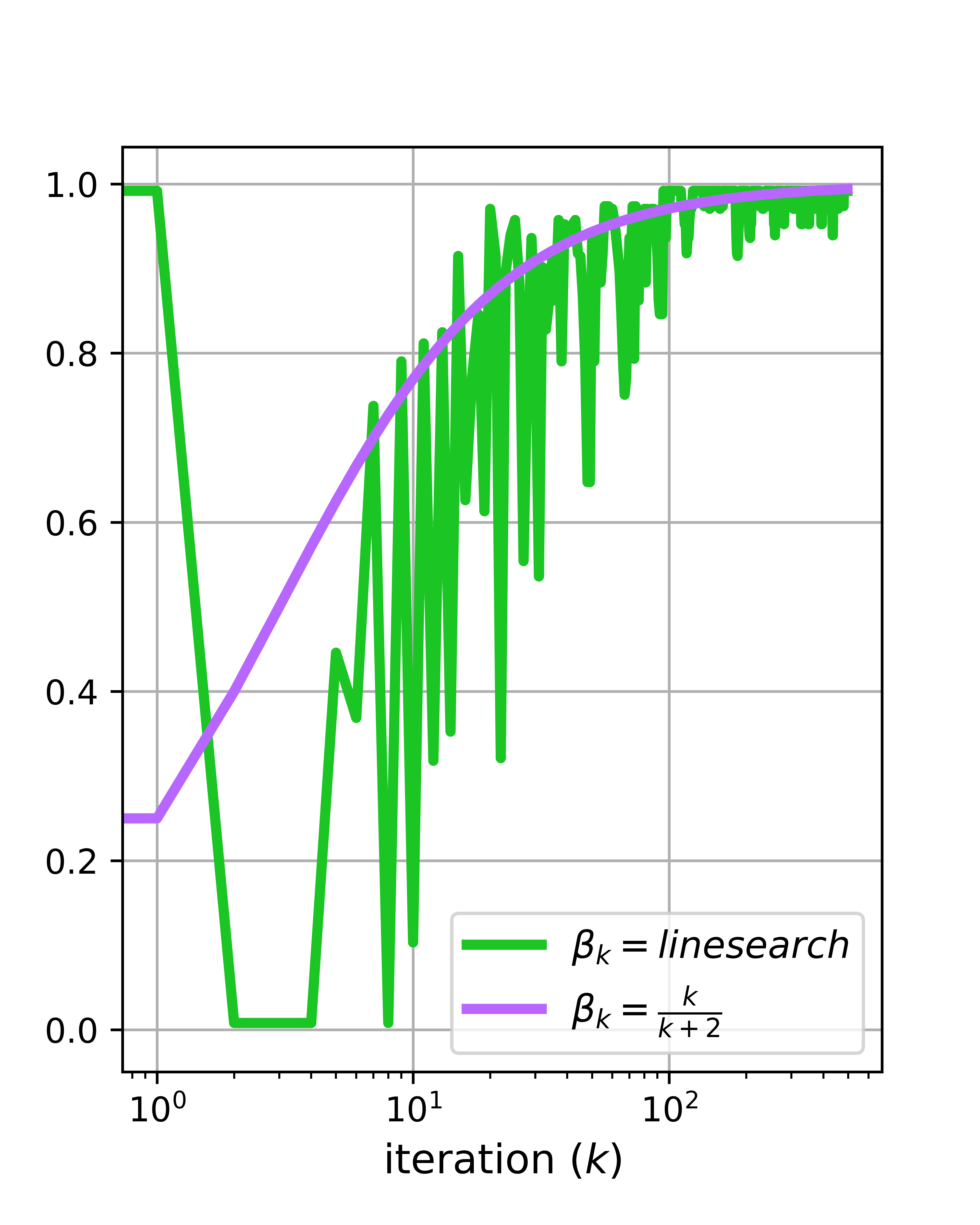}
\end{subfigure}
    \caption{\small{Maximization of the Rayleigh quotient on $M = \mathbb{S}^{d-1}$. Setting is discussed in Sec.~\ref{sec:experiments_positive}. We found that just 8 iterations of golden section search are sufficient to guarantee a steady per-iteration decrease in RAGDsDR up until a suboptimality of $10^{-9}$.}}
    \label{fig:result}
    \vspace{-4mm}
\end{center}
\end{figure}

We first consider the problem of maximizing the Rayleigh quotient $\frac{x^TAx}{2\|x\|_2^2}$ over $\R^d$, i.e. of finding the dominant eigenvector of $A\in\R^{d\times d}$. This non-convex problem can be written on the open hemisphere $\mathbb{S}^{d-1}$ (constant positive curvature) : $\argmin_{x\in\mathbb{S}^{d-1}} f(x):=-\frac{1}{2}x^TAx.$
It is well known that, in the Euclidean case, such an
objective is hard to optimize if $A$ is
high-dimensional and ill-conditioned --- and is therefore able to truly showcase the acceleration phenomenon\footnote{Indeed, high dimensional quadratics are used to construct lower bounds in~\citep{nesterov1983method}.} for convex but not necessarily strongly-convex functions, in a tight way. We choose $A = \frac{1}{d}BB^T$, where $B\in\R^{d\times n}$ has standard Gaussian entries\footnote{Inspired by PCA and linear regression, where $B$ is the design matrix ($n$ data points).}. We choose $d = 2000$ and $n=2100\approxeq d$, leading to a large condition number. In correspondence to the Euclidean case, we have $L = \lambda_{\text{max}}(A)$ and use a step-size of $1/L$ for RGD and RAGD. Also, we choose the strong-convexity modulus $\mu$~(needed parameter for RAGD) as $\lambda_{\text{min}}(A)$, again in correspondence with the Euclidean case.

\vspace{-2mm}
\paragraph{Results.} As predicted by Theorem~\ref{thm:main_theorem}, Figure~\ref{fig:result} shows that RAGDsDR is able to accelerate RGD from $\bigO(1/k)$ to $\bigO(1/k^2)$ during the first hundred iterations. The rate will eventually\footnote{This happens quite late, around iteration $k = 100$, because of the large condition number $\kappa(A) \approxeq 4000$.} become linear, due to the gradient-dominance of $f$~(Thm.~4 in~\citep{zhang2016riemannian}). In contrast, RAGD is only able to profit from acceleration \textit{at a late stage} --- before that, it is comparable to RGD.\\
We note that the choice $\beta_k=\frac{k}{k+2}$, which reduces the iteration-cost of RAGDsDR, does not influence much the empirical rate. Indeed, as shown in the figure, the geodesic search returns a result which is somehow similar. However, as also mentioned in~\citep{nesterov2018primal}, the geodesic search increases the adaptiveness of the method to curvature, providing better stability (no oscillations) and steady decrease at each iteration.  

\begin{figure}[H]
\centering
  \includegraphics[width=0.54\linewidth]{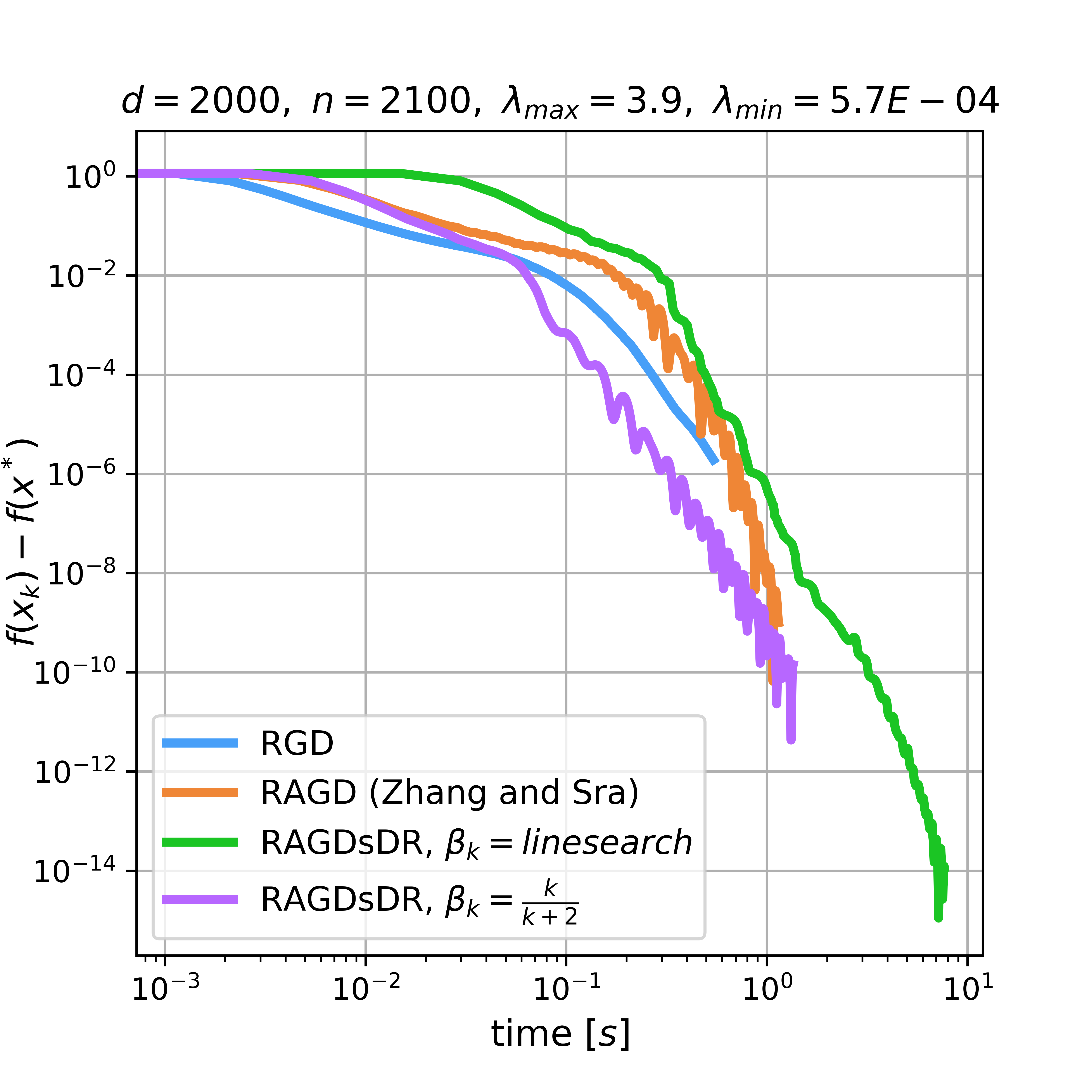}
    \caption{\small{Wall-clock time performance. Settings as Fig.~\ref{fig:result}.}}
    \label{fig:result_time}
        \vspace{-2mm}
\end{figure}

\vspace{-2mm}
\paragraph{Comment on regularization.} In the Euclidean setting, one can sometimes add a quadratic regularizer to accelerate the convergence of momentum methods designed for strongly-convex objectives. For the Rayleigh quotient problem, one may replace $A$ with $A+\gamma I_{d\times d}$, where $\gamma>0$. We note that there is typically no general principle for choosing an appropriate $\gamma$ (which is tie to generalization in machine learning). However, such a regularization technique increases the value of the strong-convexity modulus $\mu$, which speeds up optimizers designed for  strongly-convex problems. Instead, the algorithm we present in this paper provably improves over RGD in terms of gradient computations, and this effect is independent of $\mu$. To the best of our knowledge, RAGDsDR is the only Riemannian algorithm in the current literature with these features. To conclude, we also note that the derivation of accelerated rates for problems which are not strongly-convex has a long history in convex optimization (e.g. Nesterov’s 1983 seminal paper~\citep{nesterov1983method}) and arguably deserves the same attention in Riemannian optimization.
\paragraph{Comment on the wall-clock time performance.} RAGDsDR, with or without geodesic search, only requires the computation of one gradient per iteration. However, the calculation of $\beta_k$ using geodesic search (line 3 in Algorithm~\ref{alg:RAGDsDR}) increases the time complexity. The approximation of $\beta_k$ does not require additional gradients, but just a few (in this case 8 per iteration) function evaluations. For simple problems such as the ones we present in this section, the overall complexity is dominated by the call of geometric operations like the log and exponential maps (required for function evaluations along geodesics). Hence, as shown in Figure~\ref{fig:result_time}, RAGDsDR with geodesic search is de facto slower than RGD with an optimized step-size. RGD of course benefits from less geometric operations required per iteration. However, we note that (1) the practical variant of RAGDsDR is faster than RGD, and (2) for problems where the cost of a gradient computation is dominating, we would expect a significant acceleration from RAGDsDR with geodesic search. 

\subsection{Negative curvature}
\label{sec:experiments_negative}
\begin{figure}[H]
\centering
\begin{subfigure}{0.53\linewidth}
  \centering
  \includegraphics[width=\textwidth]{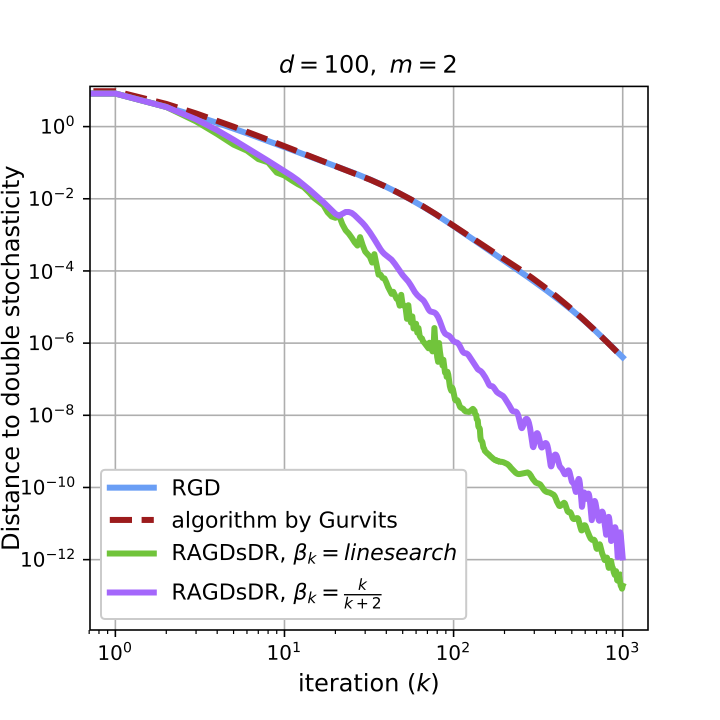}
\end{subfigure}%
\begin{subfigure}{.4\linewidth}
  \centering
  \includegraphics[width=\textwidth]{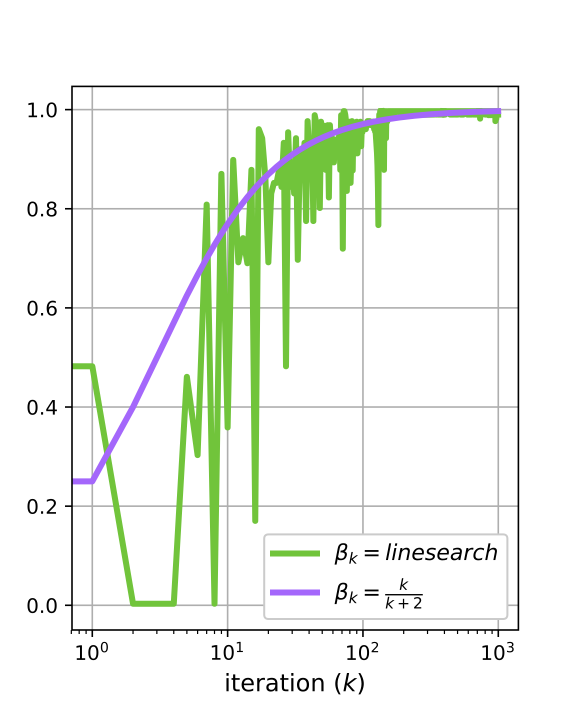}
\end{subfigure}
    \caption{\small{Scaling of a positive operator by minimizing its log-capacity. Shown is the distance to double stochasticity~(Def. 2.9 from~\citep{garg2018algorithmic}). In this metric, RAGDsDR is \textit{not necessarily a descent method}. Here we estimate $L=1$ (the smallest value that guarantees numerical stability), and note that the algorithm by~\cite{gurvits2004classical} is very similar to RGD with step $1/L$. The rate appears to be sublinear~(yet faster than $\bigO(1/k^2)$), in accordance with the complexity result in~\citep{garg2018algorithmic}.}}
    \vspace{-4mm}
    \label{fig:result_op}
\end{figure}

We study two problems on $d \times d$ symmetric positive definite matrices $\mathcal{S}^{++}(d)$. The metric $g_A(M,N) = \text{trace}(A^{-1}MA^{-1}N)$ makes $\mathcal{S}^{++}(d)$ a Riemannian manifold with negative curvature~\citep{bhatia2009positive}.
\vspace{-2mm}
\paragraph{Operator Scaling.} 
Consider an operator $T:\R^{d\times d}\to\R^{d\times d}$ defined by an $m$-tuple of $d\times d$ matrices $(A_j)_{j=1}^m$: $T(X) = \sum_{i=1}^m A_iXA_i^T$. The problem of operator scaling consists in finding $n\times n$ matrices $X$ and $Y$ such that if $\hat A_i := Y^{-1}A_iX$, then $\sum_{i=1}^m \hat A_i \hat A_i^T = \sum_{i=1}^m \hat A_i^T \hat A_i=  I_d$ (double stochasticity). Such problem is of extreme interest in theoretical computer science~\citep{garg2018algorithmic}, and has applications in algebraic complexity, invariant theory, analysis and quantum information. \cite{gurvits2004classical} showed that one can solve operator scaling by computing the \textit{capacity} of $T$, i.e. by finding $\argmin_{X\in\mathcal{S}^{++}(d)}\frac{\det(T(X))}{\det(X)}$. 
This function is non-convex in $\R^{d\times d}$, but its logarithm\footnote{$\log(\det(T(X)))-\log(\det(X))$ is geodesically convex on $\mathcal{S}^{++}(d)$. This is linked to the fact that $\log(\det(X))$ is geodesically linear (both convex and concave).} is geodesically convex on $\mathcal{S}^{++}(d)$,~\citep{vishnoi2018geodesic}.Recently,~\cite{allen2018operator} were able to exploit this property to design a competitive second-order Riemannian optimizer to solve operator scaling. Here, we instead test the performance of accelerated \textit{first-order} methods. To the best of our knowledge, there does not exist any estimate of the strong convexity constant for the log-capacity. Hence, RAGD~\citep{zhang2018towards} \textit{is not applicable} to operator scaling. Instead, we compare the performance of RAGDsDR with the algorithm by~\cite{gurvits2004classical} in Fig.~\ref{fig:result_op}, showing again a \textit{significant acceleration}.

\vspace{-2mm}
\paragraph{Karcher mean.}
Given an $n$-tuple of $d\times d$ positive definite matrices $(A_j)_{j=1}^n$, the Karcher mean is the unique positive definite solution $X$ to the equation $\sum_{i = 1}^m \log(A^{-1}_i X)=0$, where $\log$ is the matrix logarithm. This matrix average has many properties, which make its computation relevant to signal processing and medical imaging. The Karcher mean can also be written as $\argmin_{X\in\mathcal{S}^{++}(d)} f(X)=\frac{1}{2m}\sum_{i=1}^m d(A_i,X)^2$. 
\begin{figure}[H]
\centering
\begin{subfigure}{0.54\linewidth}
  \centering
  \includegraphics[width=\textwidth]{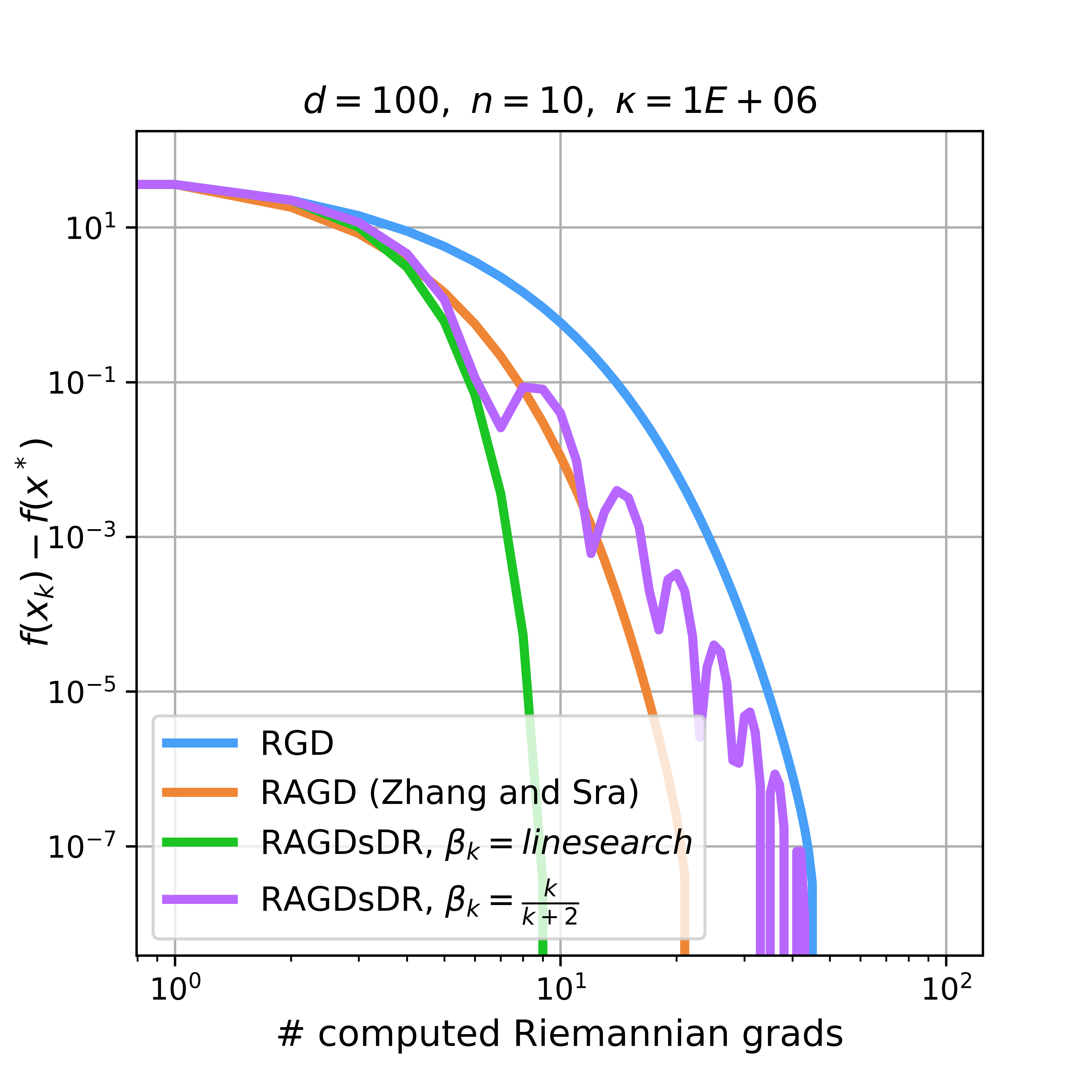}
\end{subfigure}%
\begin{subfigure}{.4\linewidth}
  \centering
  \includegraphics[width=\textwidth]{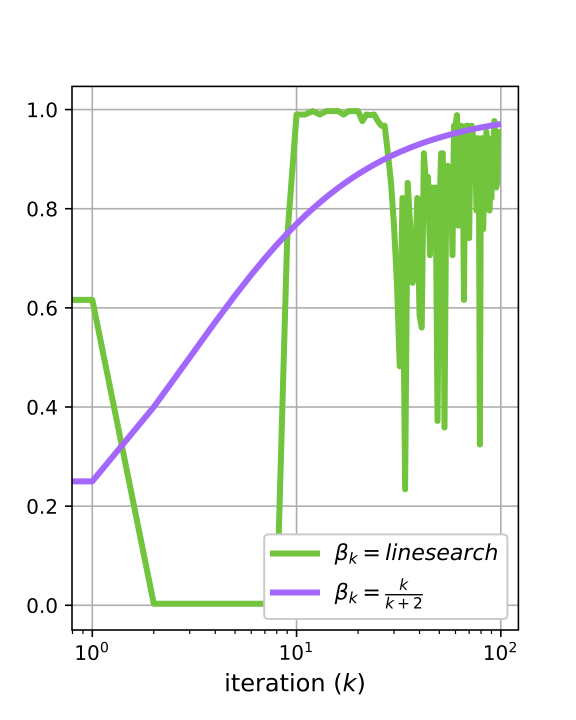}
\end{subfigure}
    \caption{\small{Performance of various optimizers on the Karcher Mean problem, as discussed in Section~\ref{sec:experiments_negative}. Performance is similar under different values for $n$ and $\kappa$. The rate appears to be linear, as predicted by~\cite{zhang2016riemannian}.}}
        \vspace{-3mm}
    \label{fig:result_ka}
\end{figure}
Clearly, $f$ is strongly-convex with modulus $\mu=1$, and $L$-smooth with modulus estimated to be around $5$~\citep{zhang2016first}. Following~\cite{zhang2016first}, we use the Matrix Mean Toolbox~\citep{bini2013computing} to generate $100$ random $100 \times 100$ positive definite matrices with fixed condition number $10^6$. In Figure~\ref{fig:result_ka}, we show that RAGDsDR (with geodesic search) is able to achieve a faster rate compared to RAGD in terms of number of iterations. Interestingly, here the choice $\beta_k=\frac{k}{k+2}$ only leads to a slight initial acceleration compared to RGD. This can be explained by looking at the values of $\beta_k$ returned by geodesic search: for the first iterations $\beta_k$ is set to a very small value --- leading to convergence in 10 iterations.

%%%%%%%%%%%%%%%%%%%%%%%%%%%
%%%%%%%%%%%%%%%%%%%%%%%%%%%%%%%%%

\section{Discussion}
\label{sec:discussion}
We proposed a novel algorithm that exploits momentum for minimizing geodesically convex and weakly-quasi-convex functions defined on a Riemannian manifold of bounded sectional curvature. We derived theoretical guarantees proving that these algorithms achieve faster rates of convergence than RGD and validated our results empirically. We conclude by contrasting our results to prior work and discussing further extensions.
\paragraph{Extension to strongly-convex case.} 
Extending our analysis to the strongly-convex case appears non-trivial. Existing analyses such as~\citep{zhang2018towards} that consider such functions, have an extra term $\frac{\mu}{2} d(y_k,x^*)^2$ in the estimate sequence, which cannot straightforwardly be dealt with in our current proof.
\paragraph{Initialization used in \citep{zhang2018towards}.}
Theorem 3 in \citep{zhang2018towards} relies on the restrictive assumption that the initialization of their algorithm is inside a ball of radius $D= \frac{1}{20 \sqrt{K}} (\frac{\mu}{L})^{\frac{3}{4}}$ centered at $x^*$. Using the strong convexity of the objective function, they are able to prove that the working domain is expanded until $\frac{1}{10 \sqrt{K}} (\frac{\mu}{L})^{\frac{1}{4}} \leq \frac{1}{10 \sqrt{K}}$. Given that we do not use strong convexity (but just convexity), this assumption would translate to a bound on the working domain of $D \leq \frac{1}{10 \sqrt{K}}$. This would in turn imply $\zeta \approx 1.003$ and $\delta \approx 0.997$. This implies that $d(M)=4 \max \lbrace \zeta-1 , 1-\delta \rbrace \approx 0.012$ and the first point of Theorem \ref{thm:flex} holds. In addition, algorithm \ref{alg:RAGDsDR} is accelerated for at least $\left[\frac{2}{0.012}\right] \approx 166$ iterations.
\paragraph{Further improvements.}
One question of practical relevance surrounds the extra error term in our rate of convergence of Theorem \ref{thm:main_theorem}. We proved that this error decays with rate $\bigO \left(d(M) /k\right)$ and that under restrictions on the working domain, our algorithm has better worst-case behaviour than RGD. However, this extra error does not allow us to claim full acceleration of our algorithm and it is a topic for future work whether such term is an artifact of our worst-case analysis. Alternatively, an interesting direction would be to study whether the extra error arises as the numerical discretization error of the ODE derived in \citep{alimisis2019continuoustime}. However, this error is practically not a significant problem since one can perform at the beginning many steps of the method with full acceleration.
\newpage
\paragraph{Acknowledgements}
The authors would like to thank professors Nicolas Boumal and Suvrit Sra for helpful discussions on the content of this paper. Gary B{\'e}cigneul
was funded by the Max Planck ETH Center for Learning
Systems during the course of this work.

\bibliography{paper} 
\bibliographystyle{plainnat}

\newpage
\onecolumn
\appendix

\begin{center}
	\textbf{\Large{Appendix: Proofs and Supplementaries}}
\end{center}
\section{Euclidean Algorithm}
\label{app:appendix A}
We restate here the Euclidean algorithm presented in \citep{nesterov2018primal}, which serves as an inspiration for developing the new Riemannian algorithm presented in this paper. One key aspect of this algorithm compared to other accelerated methods is the use of a simple 1d line search technique to obtain $\beta_k$, which makes the algorithm a descent method. This step can also be implemented efficiently in a Riemannian setting, therefore not affecting the practical aspect of the implementation of such an algorithm. The definition of $v_{k+1}$ in the following algorithm is implicit as the minimizer of $\psi_{k+1}$, while we present the same step explicitly in algorithm \ref{alg:RAGDsDR}. 

\begin{algorithm}
We recall here the euclidean algorithm, which is the basis for the Riemannian one. It is a part of algorithm 1 in \citep{nesterov2018primal}.
\caption{Accelerated Gradient Method with Small-Dimensional Relaxation
(AGMsDR)
}

\begin{algorithmic}[1]

     \STATE $A_0=0,x_0=v_0 \in \mathbb{R}^n,\psi_0(x)=\frac{1}{2} \| x-v_0 \|^2 $
     \FOR {$k \geq 0$}
     \STATE $\beta_k = \argmin_{\beta \in [0,1]} \lbrace f(v_k+\beta (x_k-v_k)) \rbrace$
     \STATE $y_k=v_k+\beta_k (x_k-v_k)$
     \STATE $x_{k+1} = y_k-\frac{1}{L} \nabla f(y_k)$
     \STATE $\frac{ a_{k+1}^2}{A_k+a_{k+1}}=\frac{1}{L}$
     \STATE $A_{k+1}=A_k+a_{k+1}$
     \STATE $\psi_{k+1}(x)=\psi_k(x)+a_{k+1}(f(y_k)+\langle \nabla f(y_k), x-y_k \rangle) $ 
     \STATE $v_{k+1}=\argmin_{x \in \mathbb{R}^n} \psi_{k+1}(x)$
     \ENDFOR
\end{algorithmic}
\end{algorithm}

%%%%%%%%%%%%%%%%%

\section{Geodesic search (equation \ref{eq:line_search_cond})}
\label{app:appendix line search}
We now examine in greater detail geodesic search in algorithm \ref{alg:RAGDsDR} (step 3) and its two main consequences summarized in equation \ref{eq:line_search_cond}.
\newline
The first condition $f(y_k) \leq f(x_k)$ follows by simply setting $\beta=1$ in the expression $f(\exp_{v_k}(\beta \log_{v_k}(x_k))$. 
\newline
For the second condition $\langle \grad f(y_k), \log_{y_k}(v_k) \rangle \geq 0$, we consider different cases depending on the value of $\beta$.
We have to take into consideration that $y_k$ is on the geodesic connecting $v_k$ with $x_k$. The derivative of the curve $\exp_{v_k}(\beta \log_{v_k}(x_k))$ with respect to $\beta$ is tangent to the geodesic and has length equal to $\| \log_{v_k}(x_k) \|$, because geodesics have constant speed. This means that the derivative at the point $y_k$ is equal to $\Gamma_{v_k}^{y_k} \log_{v_k}(x_k)$. By relying on the optimality condition of $\beta$, we distinguish the following three cases:
\setlength{\itemindent}{-0.5in}
\begin{enumerate}[label=(\roman*), leftmargin=.25in]
\itemsep0em 
    \item If $\beta_k=0$, then $\langle \grad f(y_k),\Gamma_{v_k}^{y_k} \log_{v_k}(x_k) \rangle \geq 0$ ($f(\exp_{v_k}(\beta \log_{v_k}(y_k)))$ is locally increasing on the right) and $y_k=v_k$, thus $\langle \grad f(y_k),\log_{y_k}(v_k) \rangle = 0$.
    \item If $\beta_k \in (0,1)$, then\footnote{We use Fermat's theorem for $f(\exp_{v_k}(\beta \log_{v_k}(y_k)))$.} $\langle \grad f(y_k),\Gamma_{v_k}^{y_k} \log_{v_k}(x_k) \rangle = 0$ and $\log_{v_k}(y_k)= \beta_k \log_{v_k}(x_k)$.\\
    Thus, $\langle \grad f(y_k),\frac{1}{\beta_k}\Gamma_{v_k}^{y_k} \log_{v_k}(y_k) \rangle = 0$, which implies $\langle \grad f(y_k), \log_{y_k}(v_k) \rangle = 0$.
    \item If $\beta_k=1$, then\footnote{$f(\exp_{v_k}(\beta \log_{v_k}(y_k)))$ is locally decreasing on the left.} $\langle \grad f(y_k),\Gamma_{v_k}^{y_k} \log_{v_k}(x_k) \rangle \leq 0$ and $y_k=x_k$.\\
    We deduce that $\langle \grad f(y_k),\Gamma_{v_k}^{y_k} \log_{v_k}(y_k) \rangle \leq 0$, thus $\langle \grad f(y_k),- \log_{y_k}(v_k) \rangle \leq 0$ and finally $\langle \grad f(y_k), \log_{y_k}(v_k) \rangle \geq 0$.
\end{enumerate}
In any case, the second condition is satisfied.

%%%%%%%%%%%%%%%%%%%%

\section{Proof of Lemma \ref{le:extra_error}}
\label{app:extra_error}
Consider the function $g: [0,1] \rightarrow \mathbb{R}$ defined as
\begin{equation*}
    g(t)= \langle \grad f(y_k), \Gamma_{\gamma(t)}^{y_k} \log_{\gamma(t)}(x) \rangle,
\end{equation*}
where $\gamma:[0,1] \rightarrow M$ is the geodesic connecting $y_k=\gamma(0)$ and $v_k=\gamma(1)$.
By the mean value theorem, there exists some $t_0 \in (0,1)$, such that $g(1)-g(0)=\dot g(t_0)$. This is equivalent to
\begin{align*}
\mathcal{E}_k(x) & = \langle \grad f(y_k), \log_{y_k}(x)-\Gamma_{v_k}^{y_k} \log_{v_k}(x) \rangle = \langle \grad f(y_k), \left. \frac{d}{dt} \right|_{t=t_0}-\Gamma_{\gamma(t)}^{y_k} \log_{\gamma(t)}(x) \rangle \\&= \langle \grad f(y_k), -\Gamma_{\gamma(t_0)}^{y_k} \left. \nabla_{\dot \gamma(t)} \log_{\gamma(t)}(x) \right|_{t=t_0} \rangle.
\end{align*}
The last equality holds because of a well-known property of parallel transport: 
\begin{equation*}
    \frac{d}{dt} \Gamma_{\gamma(t)}^{y_k} \log_{\gamma(t)}(x)= \Gamma_{\gamma(t)}^{y_k} \nabla_{\dot \gamma(t)} \log_{\gamma(t)}(x),
\end{equation*}
where $\nabla_{\dot \gamma}$ is the covariant derivative along $\dot \gamma$ as defined in Def.~\ref{def:covariant_derivative} (see e.g Theorem 3.3.6(vi) in \citep{robbin2011introduction}). Now we have that
\begin{align*}
\nabla_{\dot \gamma(t)} \log_{\gamma(t)}(x) & = \nabla_{\dot \gamma(t)} \left(\grad_\gamma \left(-\frac{1}{2} d(\gamma,x)^2\right)(t)\right) = \nabla_{\dot \gamma(t)} \left(\grad_\gamma \left(-\frac{1}{2} d(\gamma,x)^2 \right) \right) \dot \gamma(t) \\ &= \Hess_\gamma \left(-\frac{1}{2} d(\gamma,x^*)^2 \right) \dot \gamma(t).
\end{align*}

The derivation of the second equality can be found  in~\citep{lee2018introduction}, Chapter 11.
The last equality holds because the Hessian is by definition equal to $\nabla \grad$, and since $\gamma$ is a geodesic~\footnote{Recall that the geodesic $\gamma$, defined as $\gamma(t) = \exp(t \log_{y_k}(v_k))$, has constant velocity and the parallel transport of a tangent vector along $\gamma$ remains tangent. Thus transporting parallelly $\log_{y_k}(v_k)= \dot \gamma(0)$ from $\gamma(0)$ to $\gamma(t)$ gives the velocity at $\gamma(t)$, i.e. $\dot \gamma(t)$.}, we have
\newline
$\dot \gamma(t)=\Gamma_{y_k}^{\gamma(t)} \log_{y_k}(v_k).$
Thus
\begin{align} 
\mathcal{E}_k & (x) = \langle \grad f(y_k), \log_{y_k}(x)-\Gamma_{v_k}^{y_k} \log_{v_k}(x) \rangle \\ &= \langle \grad f(y_k),-\Gamma_{\gamma(t)}^{y_k} \Hess_\gamma \left(-\frac{1}{2} d(\gamma,x^*)^2 \right) \Gamma_{y_k}^{\gamma(t)} \log_{y_k}(v_k) \rangle
\end{align}
where we will denote the operator on the RHS by $\H := -\Gamma_{\gamma(t)}^{y_k} \Hess_\gamma (-\frac{1}{2} d(\gamma,x^*)^2) \Gamma_{y_k}^{\gamma(t)}$ (further details regarding the operator $\H$ can be found in Appendix~\ref{app:operator}).

According to Lemma 2 in \citep{alimisis2019continuoustime}, the largest eigenvalue of the operator $-\Hess_\gamma (-\frac{1}{2} d(\gamma,x^*)^2)$ is upper bounded by
\begin{equation*}
    \zeta= 
    \begin{cases}
    \sqrt{-K_{\min}} D \coth(\sqrt{-K_{\min}} D) &, K_{\min}<0 \\
    1 &, K_{\min}\geq 0
    \end{cases}
\end{equation*}
%\vspace{-1mm}
while the smallest eigenvalue is lower bounded by
\begin{equation*}
    \delta= 
    \begin{cases}
    1 &, K_{\max} \leq 0\\
    \sqrt {K_{\max}} D \cot(\sqrt{K_{\max}} D) &, K_{\max} > 0 
    \end{cases}
\end{equation*}
The eigenvalues of the operator $\H$ are exactly equal to the ones of $\Hess_{\gamma} (-\frac{1}{2} d(\gamma,x^*)^2)$, because $\Gamma_{y_k}^{\gamma(t)}=(\Gamma_{\gamma(t)}^{y_k})^{-1}$, thus the norm of the operator $\H - I_d$ satisfies
\begin{equation}
    \| \H - I_d \| \leq \max \lbrace \zeta-1,1-\delta \rbrace.
    \label{eq:bound_H}
\end{equation}
We refer the reader to the next section for the derivation of the bound on the eigenvalues of $\H$.
\newline
Now, observe that the quantity $\mathcal{E}_k(x)$ can be manipulated as follows:
\begin{align*}
& \mathcal{E}_k(x) =  \langle \grad f(y_k), \log_{y_k}(x)-\Gamma_{v_k}^{y_k} \log_{v_k}(x) -\log_{y_k}(v_k) \rangle + \langle \grad f(y_k), \log_{y_k}(v_k) \rangle \\
& \geq \langle \grad f(y_k), \log_{y_k}(x)-\Gamma_{v_k}^{y_k} \log_{v_k}(x) - \log_{y_k}(v_k) \rangle - \Tilde{\epsilon},
\end{align*}
where the last inequality holds by definition of $\Tilde{\epsilon}$ (by the geodesic search) which is such that $\langle \grad f(y_k), \log_{y_k}(v_k) \rangle \geq - \Tilde{\epsilon}$.

Using Eq.~\ref{eq:bound_H}, we finally get
\begin{align*}
    &-\langle \grad f(y_k), \log_{y_k}(x)-\Gamma_{v_k}^{y_k} \log_{v_k}(x) - \log_{y_k}(v_k) \rangle \leq \| \grad f(y_k) \| \| \H - I_d \| \| \log_{y_k}(v_k) \| \\ & \leq \| \grad f(y_k) \| \max \lbrace \zeta-1,1-\delta \rbrace D
\end{align*}
by Cauchy-Schwarz inequality.
\newline
Thus $-\mathcal{E}_k(x) \leq -\langle \grad f(y_k), \log_{y_k}(x)-\Gamma_{v_k}^{y_k} \log_{v_k}(x) - \log_{y_k}(v_k) \rangle+\Tilde{\epsilon} \leq \| \grad f(y_k) \| \max \lbrace \zeta-1,1-\delta \rbrace D+ \Tilde{\epsilon}$

%%%%%%%%%%%%%%%%%%%%%%
%%%%%%%%%%%%%%%%%%%%%%%

\section{The operator $\H$}
\label{app:operator}

An important operator in the control of the extra error arising due to the "jump" we do in our estimate sequence is $\H= -\Gamma_{\gamma(t)}^{y_k} \Hess_\gamma (-\frac{1}{2} d(\gamma,x^*)^2) \Gamma_{y_k}^{\gamma(t)} :T_{y_k} M \rightarrow T_{y_k} M$. This is actually a whole family of operators depending on $t$. Let us fix some $t$, i.e. fix one operator of the family. 
\begin{itemize}
    \item The eigenvalues of $\H$ are equal to the eigenvalues of $-\Hess_\gamma (-\frac{1}{2} d(\gamma,x^*)^2)$. Indeed, the operator $-\Hess_\gamma (-\frac{1}{2} d(\gamma,x^*)^2)$ is diagonalizable (check \citep{alimisis2019continuoustime}) and can be written as $U D U^{-1}$ in a unique way, where $D$ is diagonal formed by its eigenvalues and $U$ by its eigenvectors. Then the operator $\H$ has a unique representation in the form $\Gamma_{\gamma(t)}^{y_k} U D U^{-1} (\Gamma_{\gamma(t)}^{y_k})^{-1}=(\Gamma_{\gamma(t)}^{y_k} U) D  (\Gamma_{\gamma(t)}^{y_k} U)^{-1}$ and its eigenvalues are the diagonal entries of $D$.
    \item The largest eigenvalue of $-\Hess_\gamma (-\frac{1}{2} d(\gamma,x^*)^2)$ is less or equal than 
    \begin{equation*}
    \zeta= 
    \begin{cases}
    \sqrt{-K_{\min}} d(\gamma,x^*) \coth(\sqrt{-K_{\min}} d(\gamma,x^*)) &, K_{\min}<0 \\
    1 &, K_{\min}\geq 0
    \end{cases}.
\end{equation*}
and the smallest more or equal than
\begin{equation*}
    \delta= 
    \begin{cases}
    1 &, K_{\max} \leq 0\\
    \sqrt{K_{\max}} d(\gamma,x^*) \cot(\sqrt{K_{\max}} d(\gamma,x^*)) &, K_{\max} > 0 
    \end{cases}
\end{equation*}
Indeed, Lemma 2 in \citep{alimisis2019continuoustime} implies that
\begin{equation*}
    \delta \| \dot \gamma \| ^2 \leq\langle -\Hess_\gamma (-\frac{1}{2} d(\gamma,x^*)^2) \dot \gamma, \dot \gamma \rangle \leq \zeta \| \dot \gamma \| ^2
\end{equation*}
\textbf{for any curve} $\gamma$. Thus for a vector $v \in T_{\gamma(t)} M$ we can choose a curve $\bar \gamma$, such that $\dot {\bar \gamma} (t)=v$. This yields to the relation
\begin{equation*}
    \delta \leq \frac{\langle -\Hess_\gamma (-\frac{1}{2} d(\gamma,x^*)^2)v, v \rangle}{\| v \|^2} \leq \zeta.
\end{equation*}
By the min-max theorem, the largest eigenvalue is the maximum of $\frac{\langle -\Hess_\gamma (-\frac{1}{2} d(\gamma,x^*)^2)v, v \rangle}{\| v \|^2}$ and the smallest its minimum over all $v \in T_{\gamma(t)} M$. Thus we recover the initial estimation for the largest and smallest eigenvalue of $\H$.
\end{itemize}

%%%%%%%%%%%%%%%%%%%%%%%%%
%%%%%%%%%%%%%%%%%%%%%%%%%

\section{Proof of Theorem \ref{thm:main_theorem}}
\label{app:Main Analysis}
\begin{proof}

As in~\citep{nesterov2018primal}, the proof relies on an estimate sequence of functions, defined as
\begin{align*}
&\psi_0(x)=\frac{1}{2} \| \log_{v_0}(x) \|^2
    \\ & \psi_k(x)=\psi_k^*+\frac{1}{2} \| \log_{v_k}(x) \|^2, k\geq 1
\end{align*}
where $\psi_k^*$ is the minimum of $\psi_k$ which is yet to be specified.

The proof consists in establishing the following two inequalities -- for a suitable choice of $\psi_k^*$ -- from which one can prove the desired final result:
\begin{itemize}
    \item \textbf{C1)} $A_k f(x_k) \leq \psi_k^*$ (see definition of $A_k$ in Algorithm~\ref{alg:RAGDsDR})
    \item \textbf{C2)} $\psi_{k+1}(x) \leq \psi_k(x)+a_{k+1} (f(y_k)+\langle \grad f(y_k),\log_{y_k}(x) \rangle - \mathcal{E}_k(x))$, at least for $x=x^*$.\\
\end{itemize}

\textbf{Proof C2.}

Consider
\begin{align*}
    \psi_{k+1}^*=\psi_k^*+a_{k+1} f(y_k)-\frac{\zeta a_{k+1}^2}{2} \| \grad f(y_k) \|^2,
\end{align*}
where
\begin{equation*}
    \zeta= 
    \begin{cases}
    \sqrt{-k_{\min}} D \coth(\sqrt{-k_{\min}} D) &, k_{\min}<0 \\
    1 &, k_{\min}\geq 0.
    \end{cases}
\end{equation*}
We now have
\begin{align*}
    &\psi_k(x)+a_{k+1}(f(y_k)+\langle \grad f(y_k), \log_{y_k}(x) \rangle) \\ &=\psi_k^*+\frac{1}{2} \| \log_{v_k}(x) \|^2+a_{k+1}(f(y_k)+\langle \grad f(y_k), \log_{y_k}(x) \rangle) \\ & \geq
    \psi_k^*+a_{k+1} f(y_k)+\frac{1}{2} \| \log_{v_k}(x) \|^2+ a_{k+1} \langle \grad f(y_k),\Gamma_{v_k}^{y_k} \log_{v_k}(x) \rangle +a_{k+1} \mathcal{E}_k(x) \\ &=\psi_k^*+a_{k+1} f(y_k)+\frac{1}{2} \| \log_{v_k}(x) \|^2+ a_{k+1} \langle \Gamma_{y_k}^{v_k} \grad f(y_k), \log_{v_k}(x) \rangle+a_{k+1} \mathcal{E}_k(x) \\ & \geq \psi_k^*+a_{k+1} f(y_k)+\frac{1}{2} \| \log_{v_{k+1}}(x) \|^2- \frac{\zeta a_{k+1}^2}{2} \| \grad f(y_k) \|^2+a_{k+1} \mathcal{E}_k(x) \\ &= \psi_{k+1}^*+\frac{1}{2} \| \log_{v_{k+1}}(x) \|^2+a_{k+1} \mathcal{E}_k(x) \\
    & = \psi_{k+1}(x)+a_{k+1} \mathcal{E}_k(x),
\end{align*}
which concludes the proof of C2.\\
The last inequality follows from the definition of $v_{k+1}$ and using a trigonometric distance bound. First, we set $v_{k+1} = \exp_{v_k}(-a_{k+1} \Gamma_{y_k}^{v_k} \grad f(y_k))$ and we get
\begin{equation*}
   \log_{v_k}(v_{k+1})=-a_{k+1} \Gamma_{y_k}^{v_k} \grad f(y_k).
\end{equation*}
Thus we have
\begin{align*}
     &\frac{1}{2} \| \log_{v_k}(x) \|^2+ a_{k+1} \langle \Gamma_{y_k}^{v_k} \grad f(y_k), \log_{v_k}(x) \rangle= \frac{1}{2} \| \log_{v_k}(x) \|^2-  \langle \log_{v_k}(v_{k+1}), \log_{v_k}(x) \rangle \\ & \geq \frac{1}{2} \| \log_{v_{k+1}}(x) \|^2- \frac{\zeta}{2} \| \log_{v_k}(v_{k+1}) \|^2= \frac{1}{2} \| \log_{v_{k+1}}(x) \|^2- \frac{\zeta}{2} a_{k+1}^2 \| \Gamma_{y_k}^{v_k} \grad f(y_k) \|^2 \\&= \frac{1}{2} \| \log_{v_{k+1}}(x) \|^2- \frac{\zeta a_{k+1}^2}{2} \| \grad f(y_k) \|^2.
\end{align*}
by the basic trigonometric distance bound (lemma 5 in \citep{zhang2016first}) in the geodesic triangle $\Delta v_k v_{k+1} x$.
\newline

\textbf{Proof C1} We prove C1 by induction.
\newline
We assume that $A_k f(x_k) \leq \psi_k^*$ and we wish to prove that $A_{k+1} f(x_{k+1}) \leq \psi_{k+1}^*$.
\begin{align*}
    &\psi_{k+1}^*= \psi_k^*+a_{k+1} f(y_k)-\frac{\zeta a_{k+1}^2}{2} \| \grad f(y_k) \|^2 \\ &\geq A_k f(x_k)+a_{k+1} f(y_k)-\frac{A_{k+1}}{2 L} \| \grad f(y_k) \|^2 \\ & \geq   A_{k+1} f(y_k)-\frac{A_{k+1}}{2 L} \| \grad f(y_k) \|^2 \\
    & = A_{k+1}(f(y_k)-\frac{1}{2L} \| \grad f(y_k) \|^2) \\
    & \geq A_{k+1} f(x_{k+1}),
\end{align*}
where the last inequality follows from the definition of $x_{k+1}$ as a gradient step and $L$-smoothness of $f$.\\

\textbf{Combining C1 and C2}
Now that we have established that both C1 and C2 hold, we get
\begin{align*}
    A_k f(x_k) \leq \psi_k^* \leq \psi_k(x^*)  &\leq \sum_{i=0}^{k-1} a_{i+1}(f(y_i)+ \langle \grad f(y_i), \log_{y_i}(x^*) \rangle- \mathcal{E}_i(x))+\psi_0(x^*) \\
    &\leq \sum_{i=0}^{k-1}a_{i+1}f(x^*)+\psi_0(x^*) - \sum_{i=0}^{k-1} a_{i+1} \mathcal{E}_i(x^*)=A_k f(x^*)+\psi_0(x^*) - \sum_{i=0}^{k-1} a_{i+1} \mathcal{E}_i(x^*),
\end{align*}
where the last inequality uses the geodesic-convexity property of the function $f$.
\newline
We now have that
\begin{align*}
   &-\sum_{i=0}^{k-1} a_{i+1} \mathcal{E}_i(x^*) =  -\sum_{i=0}^{k-1} a_{i+1} \langle \grad f(y_i), \log_{y_i}(x)- \Gamma_{v_i}^{y_i} \log_{v_k}(x) \rangle= \sum_{i=0}^{k-1} a_{i+1} \langle \grad f(y_i), -\log_{y_i}(x)+\Gamma_{v_i}^{y_i} \log_{v_k}(x) \rangle \\ &\leq \sum_{i=0}^{k-1} a_{i+1} (\langle \grad f(y_i), -\log_{y_i}(x)+\Gamma_{v_i}^{y_i} \log_{v_k}(x) \rangle + \langle \grad f(y_i),\log_{y_i} (v_i) \rangle + \Tilde{\epsilon}) \\&= \sum_{i=0}^{k-1} a_{i+1} (\langle \grad f(y_i), -\log_{y_i}(x)+\Gamma_{v_i}^{y_i} \log_{v_k}(x)+\log_{y_i} (v_i) \rangle) + A_k \Tilde{\epsilon} \\ & \leq \sum_{i=0}^{k-1} a_{i+1} \|\grad f(y_i) \| \max \lbrace \zeta-1,1-\delta \rbrace D + A_k \Tilde{\epsilon} \\ &= \sum_{i=0}^{k-1} d(v_i,v_{i+1}) \max \lbrace \zeta-1,1-\delta \rbrace D + A_k \Tilde{\epsilon}  \\& \leq k \max \lbrace \zeta-1,1-\delta \rbrace D^2+ A_k \Tilde{\epsilon}
\end{align*}
The first inequality holds, because by the second property of geodesic search (equation \ref{eq:line_search_cond}). The second inequality holds by lemma \ref{le:extra_error}.
\newline
Thus we get an upper bound for the suboptimality gap:
\begin{equation}
\label{eq:bound}
    f(x_k)-f(x^*) \leq \frac{\psi_0(x^*)}{A_k}+\frac{k \max \lbrace \zeta-1,1-\delta \rbrace D^2}{A_k}+\Tilde{\epsilon}
\end{equation}

We can derive a lower bound for $A_k$ from the equation $\frac{\zeta a_{k+1}^2}{A_k+a_{k+1}}=\frac{1}{L}$ (similarly to \cite{nesterov2018primal}). Namely $A_k \geq \frac{k^2}{4\zeta L}$ and equation \ref{eq:bound} becomes
\begin{equation*}
    f(x_k)-f(x^*) \leq \frac{4 \zeta L \psi_0(x^*)}{k^2} + \frac{4 \max \lbrace \zeta-1,1-\delta \rbrace \zeta L D^2}{k}+\Tilde{\epsilon}.
    \label{eq:accelerated_rate2}
\end{equation*}

Using the fact that $\psi_0(x^*)=\frac{1}{2} d(x_0,x^*)^2$, we get:
\begin{tcolorbox}
\begin{align*}
    f(x_k)-f(x^*) \leq \frac{2 \zeta L d(x_0,x^*)^2}{k^2} + \frac{4 \max \lbrace \zeta-1,1-\delta \rbrace \zeta L D^2}{k}+\Tilde{\epsilon}
\end{align*}
\end{tcolorbox}

\end{proof}

%%%%%%%%%%%%%%%%%%%%%%%%%%%%
%%%%%%%%%%%%%%%%%%%%%%%%%%%%

\section{Proof of theorem \ref{th:convergence_thm_quasi_convex}}
\label{app:weakly-quasi-convex}
We now turn our attention to the more general class of $\alpha$-weakly-quasi-convex functions. This requires a slight modification to Algorithm~\ref{alg:RAGDsDR} by applying a restarting technique detailed in Algorithm~\ref{alg:RAGDsDR_weak}.

%%%%%%%%%%%%%%%%%%%%%%%%%%%%%%%%%%%%%%%%%%%%%%%%%%%%%%%%%%%%

%%%%%%%%%%%%%%%%%%%%%%%%%%%%%%%%%%%%%%%%%%%%%%%%%%%%%%%%%%%%

The constant $c$ in the algorithm is chosen to be bigger than $1$ ($c=2$ in \citep{nesterov2018primal}).
\begin{tcolorbox}
\begin{lemma}
\label{le:convex}
Algorithm \ref{alg:RAGDsDR} applied to an $\alpha$-weakly-convex function $f$ produces iterates $x_k$ satisfying
\begin{equation*}
    A_k (f(x_k)-f(x^*)) \leq (1-\alpha) A_k (f(x_0)-f(x^*))+\psi_0(x^*)+k \max \lbrace \zeta-1,1-\delta \rbrace D^2+ A_k \Tilde{\epsilon},
\end{equation*}
where $\psi_0(x^*)=\frac{1}{2} d(x_0,x^*)^2$.
\end{lemma}
\end{tcolorbox}
\begin{proof}
We note that both \textbf{C1} and \textbf{C2} proven in appendix \ref{app:Main Analysis} did not require convexity and we can therefore apply both inequalities to obtain:
\begin{align*}
   A_k f(x_k) \leq \psi_k^* &\leq \sum_{i=0}^{k-1} a_{i+1}((f(y_i)+\langle \grad f(y_i),\log_{y_i}(x^*) \rangle -\mathcal{E}_i(x^*)) + \psi_0(x^*) \\
   & \leq \sum_{i=0}^{k-1} a_{i+1}((1-\alpha) f(y_i)+ \alpha f(x^*) -\mathcal{E}_i(x^*)) + \psi_0(x^*) \\
   & \leq \sum_{i=0}^{k-1} a_{i+1}((1-\alpha) f(x_0)+ \alpha f(x^*) -\mathcal{E}_i(x^*)) + \psi_0(x^*) \\
   & = A_k ((1-\alpha) f(x_0)+A_k \alpha f(x^*)-\sum_{i=0}^{k-1} a_{i+1} \mathcal{E}_i(x^*) + \psi_0(x^*),
\end{align*}
where the third inequality uses the fact that the function $f$ is $\alpha$-weakly-quasi-convex.

Thus
\begin{align*}
    A_k (f(x_k)-f(x^*)) &\leq A_k (1-\alpha) (f(x_0)-f(x^*))+\psi_0(x^*)-\sum_{i=0}^{k-1} a_{i+1} \mathcal{E}_i(x^*) \\ & \leq A_k (1-\alpha) (f(x_0)-f(x^*))+\psi_0(x^*)+ k \max \lbrace \zeta-1,1-\delta \rbrace D^2+ A_k \Tilde{\epsilon}
\end{align*}
\end{proof}

\begin{tcolorbox}
\ConvergenceTheoremQuasiConvex*
\end{tcolorbox}

\begin{proof}
We first consider the first outer loop of Algorithm~\ref{alg:RAGDsDR_weak} for $i=0$. Let $\epsilon_0= f(x_0^0)-f(x^*)$. By Lemma \ref{le:convex} and the lower bound $A_k \geq \frac{k^2}{4\zeta L}$ established previously, we have that
\begin{equation*}
    f(x_k^0)-f(x^*) \leq (1-\alpha) \epsilon_0+ \frac{2 \zeta L D^2}{k^2}+d(M) \frac{\zeta L D^2}{k}+\Tilde{\epsilon}.
\end{equation*}

We want to show that the LHS is less or equal than $(1-\frac{\alpha}{c}) \epsilon_0$, therefore it suffices that
\begin{equation*}
   (1-\alpha) \epsilon_0+ \frac{2 \zeta L D^2}{k^2}+d(M) \frac{\zeta L D^2}{k}+\Tilde{\epsilon}\leq \left(1-\frac{\alpha}{c} \right) \epsilon_0.
\end{equation*}
This is equivalent to 
\begin{align*}
    &\frac{2 \zeta L D^2}{k^2}+d(M) \frac{\zeta L D^2}{k} \leq \frac{(c-1) \alpha}{c}-\Tilde{\epsilon} =:A \Longleftrightarrow \\
    & k^2-\frac{d(M) \zeta L D^2}{A} k - \frac{2 \zeta L D^2}{A} \geq 0
\end{align*}
This is satisfied if
\begin{equation*}
    k \geq \frac{d(M) \zeta L D^2}{2A}+\sqrt{\left(\frac{d(M) \zeta L D^2}{2A}\right)^2+\frac{4 \zeta L D^2}{A}}
\end{equation*}
This implies that the algorithm is first restarted after at most $N_0 = \left \lceil \frac{d(M) \zeta L D^2}{2A}+\sqrt{\left(\frac{d(M) \zeta L D^2}{2A}\right)^2+\frac{4 \zeta L D^2}{A}} \right\rceil$ iterations.

Similarly between the $i^{th}$ and the $(i+1)^{th}$ restart we have that 
\begin{equation*}
    f(x_k^i)-f(x^*) \leq (1-\alpha) \left(1-\frac{\alpha}{c} \right)^i \epsilon_0+ \frac{2 \zeta L D^2}{k^2}+d(M) \frac{\zeta L D^2}{k}+\Tilde{\epsilon} \leq \left(1-\frac{\alpha}{c}\right)^{i+1} \epsilon_0,
\end{equation*}
which is equivalent to
\begin{equation*}
\frac{2 \zeta L D^2}{k^2}+d(M) \frac{\zeta L D^2}{k} \leq \frac{(c-1) \alpha}{c} \left(1-\frac{\alpha}{c}\right)^i \epsilon_0-\Tilde{\epsilon}=:A_i,    
\end{equation*}
or
\begin{equation*}
    k \geq  \frac{d(M) \zeta L D^2}{2A_i}+\sqrt{\left(\frac{d(M) \zeta L D^2}2{A_i}\right)^2+\frac{4 \zeta L D^2}{A_i}}
\end{equation*}
Thus, between the $i^{th}$ and the $(i+1)^{th}$ restart we have at most
\begin{equation*}
    N_i= \left \lceil \frac{d(M) \zeta L D^2}{2A_i}+\sqrt{\left(\frac{d(M) \zeta L D^2}2{A_i}\right)^2+\frac{4 \zeta L D^2}{A_i}} \right\rceil \leq \left \lceil \frac{d(M) \zeta L D^2}{A_i} + \sqrt{\frac{4 \zeta L D^2}{A_i}} \right\rceil
\end{equation*}
steps ($N_i$-many steps suffice for the restart to happen).
\newline
Let $d = \log_{1-\frac{\alpha}{c}} \frac{\epsilon}{\epsilon_0}$. Then we obtain an $\epsilon$-solution using algorithm \ref{alg:RAGDsDR_weak} after $d$-many restarts. 
\newline
If algorithm 2 runs for $N$-many steps overall, we have
\begin{align*}
    N=\sum_{i=0}^d N_i \leq &\sum_{i=0}^d \left \lceil2 \frac{d(M) \zeta L D^2}{A_i} + \sqrt{\frac{4 \zeta L D^2}{A_i}} \right\rceil \nonumber \\
    &\leq d+1+ \sum_{i=0}^d \left(\frac{2 d(M) \zeta L D^2}{\frac{(c-1) \alpha}{c} \epsilon-\Tilde{\epsilon}}+  \sqrt{\frac{4 \zeta L D^2}{\frac{(c-1) \alpha}{c} \epsilon-\Tilde{\epsilon}}}\right) \left(1-\frac{\alpha}{c}\right)^{\frac{d-i}{2}}  \nonumber \\
    & = d+1+ \left(\frac{2 d(M) \zeta L D^2}{\frac{(c-1) \alpha}{c} \epsilon-\Tilde{\epsilon}}+  \sqrt{\frac{4 \zeta L D^2}{\frac{(c-1) \alpha}{c} \epsilon-\Tilde{\epsilon}}}\right) \sum_{i=0}^d  \left(1-\frac{\alpha}{c}\right)^{\frac{d-i}{2}} \nonumber \\
    & = \bigO \left( \frac{d(M)\zeta L D^2}{\alpha^2 \epsilon-\frac{c \alpha}{(c-1)}\Tilde{\epsilon}} +\sqrt{\frac{\zeta L D^2 }{\alpha^3 \epsilon-\frac{c \alpha^2}{(c-1)}\Tilde{\epsilon}}} \right)
\end{align*}
similarly to the sequence of relations at the end of Theorem 4 in \citep{nesterov2018primal}. The last equality holds because the quantity $\sum_{i=0}^d   (1-\frac{\alpha}{c})^{\frac{d-i}{2}}$ is bounded by 
a constant depending only on $\alpha$ and $c$.
\newline
Indeed 
\begin{equation*}
    \sum_{i=0}^d   \left(1-\frac{\alpha}{c}\right)^{\frac{d-i}{2}} \leq \sum_{i=-\infty}^d   \left(1-\frac{\alpha}{c}\right)^{\frac{d-i}{2}}=\sum_{i=0}^{\infty}   \left(1-\frac{\alpha}{c}\right)^{\frac{i}{2}}= \frac{1}{1-\sqrt{1-\frac{\alpha}{c}}}= \frac{1+\sqrt{1-\frac{\alpha}{c}}}{\frac{\alpha}{c}}
\end{equation*}
We conclude that
\begin{equation*}
    f(x_N)-f(x^*) \leq \epsilon \leq \bigO \left(\frac{\zeta L D^2 }{\alpha^3 N^2} \right)+d(M) \bigO\left(\frac{\zeta L D^2}{\alpha^2 N} \right)  +\frac{c}{(c-1)\alpha}\Tilde{\epsilon}.
\end{equation*}
\end{proof}

\end{document}